%% file: nb.tex
\documentclass[12pt,a4paper]{amsart}
\usepackage{amsfonts,amssymb,amsmath}
\usepackage[all]{xy}
\usepackage[dvips,a4paper]{geometry}
\usepackage{graphicx}
\usepackage{longtable}
\theoremstyle{plain}
\newtheorem{lemma}{Lemma}
\newtheorem{theo}[lemma]{Theorem}
\newtheorem{propo}[lemma]{Proposition}

\theoremstyle{definition}

\theoremstyle{remark}
\newtheorem{remark}{Remark}
\newdir{ >}{!/8pt/@{ }*@{>}}
\newdir{< }{!/-8pt/@{ }*@{<}}

\newcommand{\C}{\mathbb{C}}
\newcommand{\R}{\mathbb{R}}
\newcommand{\CP}{\mathbb{CP}}

\newcommand{\RP}{\mathbb{RP}}

\newcommand{\CC}{\mathrm{CC}}
\newcommand{\fix}{\mathrm{Fix}}
\newcommand{\fbar}{\overline{f}}
\newcommand{\from}{\colon}
\newcommand{\Chi}{\mathcal{X}}
\newcommand{\minus}{\smallsetminus}
\newcommand{\crit}{\mathrm{Crit}}
\newcommand{\barCC}{\overline{\mathrm{CC}}}
\begin{document}

\pagenumbering{arabic}

\title[Central configurations, symmetries and fixed points]{%
Central configurations, symmetries and fixed points
}

\author{Davide L.~Ferrario}

\address{
Dipartimento di Matematica del
 Politecnico di Milano\\
Piazza Leonardo da Vinci 32 --
 20133 Milano (I) \\
Current address: Max-Planck-Institut f\"ur Mathematik\\
Vivatsgasse, 7 - 53111 Bonn (DE)
 }
\email{ferrario@mate.polimi.it}

\date{\today}

\begin{abstract}
Planar central configurations can be seen as critical points of 
the reduced potential or solutions of 
a system of equations. 
By the homogeneity and invariance of the potential with respect to $SO(2)$,
it is possible to see that the $SO(2)$-orbits of central configurations 
are fixed points of a suitable map $f$.
The purpose of the paper is to define this map and to derive 
some properties using topological fixed point theory.
The generalized Moulton-Smale theorem  for collinear 
configurations is proved, together with  some estimates on the number 
of central configurations in the case of $3$ bodies,
using fixed point indexes.  Well-known results such as the 
compactness of the set of central configuration can also be 
proved in an easy way in this topological framework.
At the end of the paper some tables of 
(numerical) planar central configurations of $n$ 
equal masses with Newtonian potential are given,
for $n=3,\dots, 10$. They have been computed as the fixed points
of a suitable self-map of $\R^{2(n-2)}$.
\end{abstract}

 \subjclass{
 Primary 
 70F10 
 Secondary 
37C25 
 }

 \keywords{%
Central configurations, fixed points, fixed point index,
Moulton configurations
}

\maketitle


\section{Planar central configurations}
Consider the complex plane $\C$ and $n$ real positive numbers
$m_1$ , $m_2$ , $ \dots$ , $m_n$, with $n\geq 3$. We denote by
$\Chi$ the hyperplane in $\C^n$ defined by the equation
$\sum_i m_iz_i=0$, where $z=(z_1,z_2,\dots,z_n)$ denotes
a point in $\C^n$. Let $\hat \Chi$ be the open subset of 
$\Chi$ defined by the condition $i\neq j \implies z_i\neq z_j$.
That is, if $\Delta_{ij}$ denotes the subspace of 
equation $z_i=z_j$ and $\Delta=\cup_{i<j} \Delta_{ij}$ the 
collision set, it is $\hat \Chi = \Chi \minus \Delta$.
It is the configuration space of $n$ point particles in $\C$ and
center of mass in $0$. An element in $\hat \Chi$ is called
a configuration.
Consider some real coefficients $m_{ij}$, for $i\neq j$, $i$, 
$j \in \{1,\dots n\}$  and a regular function $\phi\from \R_+
\to \R$ defined on the positive semiline of $\R$.
Without loss of generality we can assume $m_{ij}=m_{ji}$ for 
all $j\neq i$.
Let $\alpha$ be a real number and let $U\from \hat \Chi \to \R$
be a potential function 
be of type
\begin{equation*}
U = \sum_{i<j} m_{ij} \phi(z_i-z_j),
\end{equation*}
where $\phi(z) = |z|^{\alpha+2}$ or $\phi(z)=\log|z|$ (in case $\alpha=-2$).
For example if 
$m_{ij} = m_i m_j$ and $\alpha=-3$ then this is 
the potential of the Newtonian $n$-body problem with $\phi(z)=|z|^{-1}$
or the Thompson vortex problem if $\alpha=-2$ and then
 $\phi(z) = \log|z|$.
The charged $n$-body problem is obtained by
setting $\alpha=-3$ and $m_{ij} = m_im_j - q_iq_j$ where $q_i$ 
are the electrostatic charges of the masses.

The $n$-body problem concerns the motion of $n$ particles 
of masses $m_i$ and potential $U$; the Newton equations are 
therefore $m_i\ddot z_i = \frac{\partial U}{\partial z_i}$.
A \emph{central configuration} is a configuration $z\in \hat \Chi$
such that there exists a non-zero real scalar $k$ 
such that  for $i=1,\dots n$
\begin{equation*}
k m_i z_i = \frac{\partial U}{\partial z_i}. 
\end{equation*}
Central configurations can be seen as critical points of  
the restriction of $U$ to the ellipsoid in $\hat \Chi$ 
of equations $\sum_i m_i|z_i|^2 = 1$.
Furthermore, they play an important role in the theory
of periodic orbits in $n$-body problems because they
yield homographic solution and are topological bifurcations 
of the energy and angular momentum level sets in $\hat\Chi$.
Details and further reading on the topic can be found for example in
\cite{Mo90,PSSY,MC96,Sm70I}.

Let $\hat f\from \hat \Chi \to \Chi$ be the map defined by
$\hat f(z) = w \in \Chi$, where
\begin{equation*}
w_i = \sum_{j\neq i} m_i^{-1} m_{ij} |z_i - z_j|^{\alpha}(z_i-z_j).
\end{equation*}
It is not difficult to see 
$\dfrac{\partial U}{\partial z_i} = (\alpha+2) m_i w_i$, 
so that a configuration $z\in \hat\Chi$ is central
if and only if there is a real number $\lambda\neq 0$
such that 
$\hat f(z) = \lambda z$,
where $\lambda = (\alpha+2)^{-1}k$.
Let $\CC$ denote the space of all the 
central configurations in $\hat \Chi$ and 
$\CC_1$ its intersection with the ellipsoid $E \subset \Chi$ of equation
$\sum_i m_i |z_i|^2=1$.
By the $O(2)$-invariance of $U$, $\CC_1$ is $O(2)$-invariant
in $\hat \Chi$, where $O(2)$ is the orthogonal group
of the plane $\C$. We are interested in the quotient spaces 
$\CC_1/O(2)$ and $\barCC = \CC_1/SO(2)$. 

Consider the projection $p\from \C^n\minus\{0\} \to \CP^{n-1}$ onto the complex
projective space. A point of $\CP^{n-1}$ of homogeneous
coordinates $z_i$ is denoted by $[z]=[z_1:z_2:\dots:z_n]$.
Then $p$ projects $\Chi\minus \{0\}$ onto the hyperplane in $\CP^{n-1}$ 
of equation 
$\sum_i m_i z_i$, which we denote simply by $\CP^{n-2}$.
Let $X$ be the image of $\hat \Chi$ in $\CP^{n-2}$.
With an abuse of terminology we will consider maps
defined on open dense subsets of their domains. Such a subset
of points in which the map is properly defined will 
be clear from the context. 
The map $\hat f$ induces a map
$f\from X \to \CP^{n-2}$. It is clear that a central
configuration projects to a fixed point of $f$.

\begin{lemma}\label{lemma:2}
The projection $p$ induces 
a homeomorphism $p\from \barCC = \CC_1/SO(2) \approx \fix(f)$.
\end{lemma}
\begin{proof}
The map $p$ induces a 
continuous map $\CC_1 \to \fix(f)$, 
and a injective continuous 
map $\bar p\from \CC_1/SO(2) \to \fix(f)$.
To see that it is surjective, consider a configuration
$z$ such that $[z] = f[z]$, i.e. a configuration
$z=[z_1:z_2:\dots:z_n]$ such that there exists 
$\lambda \in \C^*$ and for every $i=1\dots n$
\begin{equation*}
\lambda z_i =
\sum_{j\neq i}
m_i^{-1} m_{ij} |z_i-z_j|^{\alpha}(z_i-z_j).
\end{equation*}
This implies 
\begin{equation*}
\lambda m_i |z_i|^2 =
\sum_{j\neq i}
m_{ij} |z_i-z_j|^{\alpha}(|z_i|^2-z_j\bar z_i) 
\end{equation*}
\begin{equation*}
\implies
\lambda \sum_i m_i |z_i|^2 =
\sum_{j\neq i}
m_{ij} |z_i-z_j|^{2+\alpha},
\end{equation*}
and therefore $\lambda$ is real, i.e. $[z]$ is a projection
of a central configuration.
Moreover, $\bar p$ is a closed map and hence a homeomorphism.
\end{proof}

Let $I=\sum_i m_i |z_i|^2$ the inertia of $z$. 
In case $\alpha\neq -2$ (not the logarithmic case), 
the proof of the 
previous lemma implies that 
\begin{equation}\label{eq:1}
\lambda  =  2\dfrac{U}{I}
\end{equation}
and hence, if $m_{ij}>0$ for all $i$, $j$ 
that $\lambda>0$. If $\alpha=-2$ the the same holds,
since 
\begin{equation}\label{eq:2}
\lambda = \dfrac{\sum_{j\neq i} m_{ij} }{I}.
\end{equation}

\begin{lemma}\label{lemma:shub}
If for every $i\neq j$ the coefficient $m_{ij}$ is not zero, 
and $\alpha<-1$, 
then $\CC_1$ is compact. 
\end{lemma}
\begin{proof}
Because $p\from E \to \CP^{n-2}$ is a proper map, 
it suffices to show that $\barCC$ is compact,
hence that $\fix(f)$ is compact.
First, let $X_2\subset \CP^{n-2}$ be the subset of 
$\CP^{n-2}$ consisting of all the points 
such that there is at most one pair of 
particles colliding $z_i=z_j$.
The map $f$ can be extended continuously to $X_2$:
If $z$ tends to a point $z'$ in $X_2\minus X$ with a collision
in $i$ and $j$, then 
the image of $z$ tends to 
the point $[w]=[w_1:\dots: w_n]$ with 
$w_i =m_jm_{ij} $ and $w_j=-m_im_{ij} $ and otherwise $0$.
Let $W_{ij}$ denote such point.
Provided that $m_{ij}\neq 0$ 
such a point $W_{ij}$ is in $\CP^{n-2}$ and it is different from
$z'$ (we assume that the masses $m_i$ are positive 
and hence it cannot be that $m_i =- m_j$):
\begin{equation*}
W_{ij} = [ 0:\dots: 0: \underset{i}{m_j} 
: 0 : \dots : 0 : \underset{j}{-m_i} :
0 :\dots : 0 ].
\end{equation*}
Now suppose that $[z]$ tends to a multiple collision $[c]$ in
$\CP^{n-2}\minus X_2$. Let $\Gamma$ be the set of indexes
$(i,j)$ such that the $i$-th particle collides with 
the $j$-th particle in $[c]$. It is not difficult to show
as above that $f(z)$ tends to a subset of the
projective subspace $\hat c \subset 
\CP^{n-2}$ spanned by the points
$W_{ij}$ with $(i,j) \in \Gamma$.
But  $\hat c$ is a closed subspace which does not contain
$[c]$, hence by continuity there is a neighborhood
of $[c]$ in $\CP^{n-2}$ without fixed points of $f$.
Therefore the collision set is contained in a 
fixed point free neighborhood, and hence $\barCC$ is closed
in $\CP^{n-2}$. Being a closed in a compact, it is compact
(compare with the proof of Shub, for the Newtonian case \cite{shub},
and with the estimates of Buck \cite{Bu90}).
\end{proof}

\begin{lemma}\label{crit}
Let $u\from X \to \R$ be the map defined by
$u[z_1:\dots:z_n] =  U I^{-1-\alpha/2}$ if $\alpha\neq -2$,
or $u[z_1:\dots:z_n] = e^U I^{-1/2\sum_{i<j} m_{ij}}$ if $\alpha=-2$.
Then $\barCC = \crit(u)$, where $\crit(u)$ denotes 
the set of critical points of $u$.
\end{lemma}
\begin{proof}
In both cases $u$ is well-defined 
on $X$.
Furthermore, $\CC_1$ is the set of critical points of $U$ restricted
to $E$; 
hence it is the set of critical points of $u\pi$,
where $\pi\from E \to X$ is the projection, restricted to $E$.
The projection $\pi$ is a submersion,
hence $\barCC = \crit(u)$.
\end{proof}

Let $C$ be the group of order $2$ acting on $\CP^{n-2}$ by
conjugation on coordinates. The space fixed by $C$ is the 
space of collinear configurations, and it is homeomorphic
to $\RP^{n-2}$. Let $X^C$ be its intersection with $X$.
The map $f$ is equivariant with respect to the action of $C$,
therefore it induces
a map $f^C\from X^C \to \RP^{n-2}$; its fixed points,
by lemma \ref{lemma:2} are the collinear central configurations.

We list some  known results. Some of them were proved  in the Newtonian case
($m_{ij} = m_i m_j$ and $\alpha=-3$), but the techniques worked 
in the same way in the general case.
We understand that 
central configurations are counted in $\CC_1/O(2)$:
For every $n\geq 3$ there are exactly 
$n!/2$ collinear configurations (Moulton; see Smale \cite{Sm70II}
for a proof using critical point theory).
If $n=3$ then there is just one 
non-collinear central configuration. 
If the masses are equal,
then there are exactly $19$ non-collinear central configurations for $n=4$ 
(Albouy \cite{A95}). 
For every $n$ there are at least $n-2$ non-collinear central
configurations (McCord \cite{MC96}).
If the potential $U$ is a Morse function, there are at least $n!h(n)/2$
central configurations, where $h(n)=\sum_{i=3}^n 1/i$ (McCord \cite{MC96}).
The Euler characteristic $\chi(X)$ is $(-1)^n(n-2)!$, 
hence by Morse theory in this case the alternating sum 
$\sum_k (-1)^k \nu_k = (-1)^n(n-2)!$, 
where $\nu_k$ denotes the number of critical points of index $k$.

\section{Collinear configurations}

We have seen that the collinear central configurations
are the fixed points of $f^C\from X^C \to \RP^{n-2}$.
We show a proof of Moulton theorem using fixed point theory,
instead of critical point theory. In some sense it is closer to the original
proof of Moulton.

\begin{propo}\label{propo:index}
If $m_{ij}>0$ for every $i$, $j$ and $\alpha<-1$ then
every fixed point of $f$ is isolated and its fixed point index 
is $1$.
\end{propo}
\begin{proof}
Let $E^C\subset \Chi^C$ be the ellipsoid  of equation 
$\sum_{i}m_iz_i$. Because $E^C \to \RP^{n-2}$ is a covering map
and $\lambda > 0$ (by equations \ref{eq:1} and \ref{eq:2}),
the map  $f'\from  E^C\minus \Delta \to E^C$ defined by 
$w \mapsto p$, with $p_i = \frac{w_i}{\sqrt{I}}$,
where as above $I=\sum_j m_jw_j^2$ and $w_j$ is defined in equation
\ref{eq:defi}, 
is a lifting of $f$.
Hence 
the Jacobian of $f$ at a point $[c]\in X^C$
is the same as the jacobian of $f'$ at a pre-image of $[c]$ in $E^C$.
We are going to show that if $x\in E^C$ is a central configuration
(that is, a fixed point for $f'$), then for every vector
$v$ of the tangent space of $E^C$ in $x$ 
(endowed
with the kinetic scalar product: $v\cdot v' = \sum_i m_i v_i v'_i$)
the inequality 
$v\cdot D(f')v < 0$ holds, where $D(f')$ denotes the 
differential of $f'$ at $x$. 
From this the claim follows, since all the eigenvalues of $D(f')$ 
need be negative.
The map
$f'$ is the  composition of $\hat f$ and the projection $p$,
hence $D(f') = D(p) D(\hat f)$.
The derivatives of $w_i$ are
\begin{equation}
\dfrac{\partial w_i}{\partial z_k} =
\begin{cases}
(\alpha+1) \sum_{j\neq i} m_i^{-1}m_{ij} |z_i - z_j|^{\alpha} 
& \mbox{if $i=k$}\\
-(\alpha+1) m_i^{-1} m_{ik} |z_i - z_k|^{\alpha}
&
\mbox{if $i\neq k$}.
\end{cases}
\end{equation}
Hence
\begin{equation}\label{eq:vivj}
\sum_{k=1}^n 
\dfrac{\partial w_i}{\partial z_k} 
=
(\alpha + 1) 
\sum_{j\neq i} m_{1}^{-1} m_{ij} |z_i - z_j|^{\alpha}(v_i - v_j).
\end{equation}
The derivatives of $p$ are
\begin{equation}
\dfrac{\partial p_i}{\partial w_k} 
=
\begin{cases}
I^{-3/2} \sum_{j\neq i} m_j w_j^2
&
\mbox{if $i=k$}\\
-I^{-3/2} m_kw_iw_k 
&
\mbox{ if $i\neq k$}\\
\end{cases}
\end{equation}
Hence 
\begin{equation}
\sum_{i=1}^n
m_i v_i 
\dfrac{\partial p_i}{\partial w_k}
=
m_k v_k \dfrac{\partial p_k}{\partial w_k}
+
\sum_{i\neq k}
m_i v_i 
\dfrac{\partial p_i}{\partial w_k}.
\end{equation}
Now, $v$ belongs to the tangent space at $x$ and if 
$x$ is a central configuration then
$\sum_i m_i w_i v_i = 0$, 
therefore
\begin{equation}
\sum_{j\neq k}
m_i v_i w_i = - m_k v_k w_k,
\end{equation}
hence 
\begin{equation}
\sum_{i\neq k}
m_i v_i 
\dfrac{\partial p_i}{\partial w_k} =
-I^{-3/2}\sum_{i\neq k}
m_i v_i w_i m_k w_k =
-I^{-3/2} m_k w_k.
\end{equation}
Therefore 
\begin{equation}
\sum_{i=1}^n
m_i v_i 
\dfrac{\partial p_i}{\partial w_k}
= m_k v_k I^{-1/2}.
\end{equation}
Now, this implies that the claim is true if and only if 
for every $v$ in the tangent space 
of $x$ the inequality 
\begin{equation}
\sum_{i,k} m_i v_i \dfrac{\partial w_i}{\partial z_k} v_k <0
\end{equation}
is true.
But by equation 
\ref{eq:vivj} 
the latter is equal to 
\begin{equation}
(\alpha + 1) 
\sum_{j\neq i} m_{ij} |z_i - z_j|^{\alpha}(v_i - v_j)v_i
=
(\alpha+1) \sum_{i<j} m_{ij} |z_i - z_j|^{\alpha}(v_i - v_j)^2,
\end{equation}
which is negative because by assumption $\alpha+1<0$  and $m_{ij}>0$.
This finishes the proof.
\end{proof}

\begin{propo}\label{propo:moulton}
If $m_{ij}>0$ for every $i<j$ and if $\alpha<-1$ then
the fixed point index of $f^C$ 
in each of the $n!/2$ components of $X^C$ is $1$.
\end{propo}
\begin{proof}
The space $X^C$ has $n!/2$ connected components
(see for example Smale \cite{Sm70II}).
Consider, for every 
$t\in I$, the map defined by $f_t(z)=w$, with
\begin{equation}\label{eq:defi}
w_i = \sum_{j\neq i} m_i^{-1}(tm_{ij} + (1 - t)m_im_j) |z_i - z_j|^{t\alpha -2 + 2t}(z_i - z_j)
\end{equation}
It is the self-map corresponding to the collinear 
$n$-body problem
with parameters $\alpha'=t\alpha -2 + 2t$ and 
$m_{ij'} = tm_{ij} + (1 - t)m_im_j$ and masses $m_i$.
By lemma \ref{lemma:shub}, 
$\fix(f_t)$ is compact for every $t\in I$.
Therefore 
$f_t$ yields  a compactly fixed homotopy from $f_1=f$
to a map $f_0\from X^C \to \RP^{n-2}$.
The map $f_0$ is the self-map arising from the logarithmic potential
collinear problem with masses $m_i$.
Now we are going to define a compactly fixed homotopy $g_t$ 
from $f_0=g_0$ to a map $g_1$ 
corresponding to the self-map of the problem with all the masses equal to $1$
and the last mass equal to $m$.
Let $m'_i = (1-t)m_i + t$ and $m'_n=(1-t)m_i + tm$; and let the hyperplane 
$\RP^{n-2}_t$ be defined by the equation $\sum_i m'_i z_i=0$
in $\RP^{n-1}$ with homogeneous real coordinates $z_i$.
By definition $X\subset \RP^{n-2} = \RP^{n-2} = \RP^{n-2}_0$.
Then there is a family of homeomorphisms
$\varphi_t\from \RP^{n-2}_0 \to \RP^{n-2}_t$.
Let $G_t\from \RP^{n-2}_t \minus \bar \Delta \to \RP^{n-2}_t$
be the self-map yielded by the $n$-body problem with masses
$m_{i}'$, $m_{ij}=m_im_j$ and $\alpha=-2$.
The composition $g_t = \varphi_t G_t \varphi_t^{-1}$ 
is a compactly fixed homotopy from $g_0=f_0$ 
to a self-map conjugated to the self-map  $h_m$ of the 
problem with logarithmic potential and masses $1$, $1$, \dots, $m$.
Now we prove, somehow as in the proof of Moulton \cite{moulton}, 
that the fixed point index of 
$h_m$ in each component of $X^C$ is $1$
for every $m$.
By induction on the number of bodies $n$. If $n=3$, then this is true.
Otherwise, consider the map $h_m$. 
It is possible to let $m$ be equal to zero (the case
of the infinitesimal mass), and everything that we have 
done in the previous section can be carried out literally.
By looking at the   
proof of lemma \ref{lemma:shub}, it is easy to see that 
if $m=0$ then $h_0$ is compactly fixed, like in the case $m>0$.
Therefore
$h_m$ is compactly fixed homotopic to $h_0$, 
and 
the fixed point indexes of 
 $h_0$ coincide with the fixed point indexes of
$h_m$.
So consider the projection $\pi\from \RP^{n-2} \to \RP^{n-3}$
which sends $[z_1:\dots : z_n]$ to $[z_1:\dots:z_{n-1}]$.
It is well-defined, because if $m_n=0$ 
then $\sum_i m_i z_i = 0 \implies m_1+\dots+m_{n-1}=0$.
The map $h_0$ induces a map $h'$ arising from the $(n-1)$-body problem
with masses $1$, \dots , $1$ like in the following diagram.
\begin{equation*}\xymatrix{%
X \ar[r]^{h_0} \ar[d]^{\pi} & \RP^{n-2} \ar[d] \\ 
X' \ar[r]^{h'} & \RP^{n-3}\\
}
\end{equation*}
By the induction hypothesis the fixed point indexes of 
$h'$ in the $(n-1)!/2$ components of $X$' 
are $1$. There are a finite number of fixed point,
being isolated and contained in a compact.
Moreover $\pi\from X\to X'$ is a fiber bundle with fibers 
equal to $\R$ minus $n-1$ points (the collisions). 
Consider a fixed point $[z_1:\dots,z_{n-1}]$ of $h'$.
Let $F\approx \R$ be its pre-image under $\pi$.
By equation \ref{eq:2}, we have $w_i = \lambda z_i$,
with $\lambda = n(n-1)/I$.
Hence the induced map 
$h_0|F$ is 
\begin{equation}
x\in \R \mapsto \dfrac{I}{n-1} \sum_{j=1}^{n-1} (x - z_j)^{-1}.
\end{equation}
It is easy to see that $h_0|F$ has just one fixed point of index $1$
in each of the $n$ components of $F$.
By the product formula for the fixed point index 
\cite{Br71}, this means that the fixed point index of $h_0$ in each of the 
$n (n-1)!/2 = n!/2$  components of $X$  is
$1$.
This completes the proof.
\end{proof}

\begin{theo}\label{propo:genmoulton}
If $\alpha<-1$ and the coefficients $m_{ij}$ are non-zero and have 
the same sign for every $i$, $j$, then
there are exactly $n!/2$ collinear central configurations.
\end{theo}
\begin{proof}
Assume $m_{ij}>0$ for every $i$, $j$.
Proposition \ref{propo:moulton} implies that 
the fixed point index of $f$ 
in every component $V$ 
of $X$ is $1$. By proposition \ref{propo:index}
and \ref{lemma:shub}, 
in each component 
$V$ there are a finite number of fixed points of index $1$.
By the additivity property of the fixed point 
index, this implies that there is exactly
one fixed point in each component, henceforth that there
are exactly $n!/2$ collinear central configurations.
If $m_{ij} < 0 $ for every $i$, $j$, then apply the same 
argument to the map $-f$ obtained by taking $-U$ instead of $U$.
Since $-f=f$, $\fix(f) = \fix(-f)$ hence the claim.
\end{proof}

\section{Three bodies}
The central configurations in the Newtonian $3$-body problem 
are the three Euler collinear configurations and the equilateral 
Lagrange configuration. More generally, 
the case of three charged
bodies has been done in \cite{PSSY} in case $\alpha=-3$.
Here we give some bounds on the number of solutions, 
assuming as above that 
$\alpha<-1$ (that is,
the collisions are singularities for the field). 
The space $X$ is homeomorphic to $\CP^1 = S^2$ minus three
points (the three double collisions). 
The map $f$
can be extended to $S^2$, since there are no triple collisions
(apply the same argument in the proof of lemma \ref{lemma:shub}).
The action of the conjugation group $C$ 
yields the reflection along the equator of $S^2$.
A configuration is non-collinear if and only if it is not
fixed by $C$. 

\subsection{Collinear solutions}
By the generalization of 
Moulton theorem \ref{propo:genmoulton} if the coefficients 
$m_{ij}$ have the same sign 
for every $i$, $j$ then there are exactly $3$ collinear central 
configurations. 
Otherwise, 
assume that some $m_{ij}$ are positive and some negative 
and none is zero.
Up to rearranging
indexes and changing sign, we can suppose $m_{12}>0$, $m_{13}>0$ and 
$m_{23}<0$,
thus that $f^C$ is compactly fixed 
homotopic to the map $\varphi$ corresponding to the problem 
with 
$m_1=m_2=m_3=1$, $m_{12}=1$, $m_{13}=1$ and $m_{23}=-1$  and 
$\alpha=-2$.
In real projective coordinates, we have 
\[
\varphi[z_1:z_2:z_3] =
[3z_1(z_2 - z_3):-(z_1-z_3)^2:(z_1-z_2)^2];
\]
in the affine chart $[t:1:-1-t]$ the map $\varphi$ 
can be written as a map $\R \to \R$ 
\[
t \mapsto -3\dfrac{t(2+t)}{(2t+1)^2},
\]
which has degree $0$.
Thus the degree of $f^C$ is equal to $0$.
This means that there is always at least one
collinear central configuration. 
In this case, opposite to the case $m_{ij}>0$, 
the number of fixed points might be greater than  $1$ (the fixed 
point index of $f^C$).
According to \cite{PSSY}, for $\alpha=-3$ the fixed points
can be any number from 1 to 5, depending on the coefficients.

Now consider the case in which some $m_{ij}$ vanish.
If all the $m_{ij}$ vanish, then the map $f$ is even not defined
and the problem is totally degenerate.
If two of the $m_{ij}$'s vanish, then the map $f$ is equal 
to a constant, so that there is exactly one fixed point.
In this case the three body problem reduces to a Kepler problem 
with a third mass which does not interact with the first two
and there is not much physical relevance.
So assume that 
only one of the $m_{ij}$ vanishes: without loss of generality
it is 
$m_{23}$. The homotopies
yielded by the variation of the coefficients 
now need not to be compactly fixed away
from the collisions, but are defined also on 
the collision set.
There are two cases: either $m_{12}$ and $m_{13}$ have the 
same sign, or not.
If they have the same sign,
then the map $f^C$ up to sign is homotopic to the map
$f_3^C$ obtained by setting $m_i=1$ and $m_{12}=m_{13}=1$ and 
$\alpha=-2$.
In coordinates, it is 
\[
f_3=[3z_1:z_3-z_1:z_2-z_1],
\]
which has affine form
\[
t \mapsto -3\dfrac{t}{2t+1}.
\]
This map has degree $-1$, therefore there are at least
$2$ fixed points. 
The only collision fixed by $f_3$ is 
the collision of the particle $2$ with
the particle $3$, 
therefore there is always $1$ fixed point
of $f^C$ 
corresponding to a central configuration.
In case $m_{12}$ and $m_{13}$ have different signs,
then
$f^C$ is homotopic to the map $f^C_4$ obtained by setting
$m_i=1$ and $m_{12}= - m_{13} = 1$ and $\alpha=-2$. 
In coordinates, it is 
\[
f_4=[z_2-z_3:z_3-z_1:z_1-z_2].
\]
This map has degree
$1$, and does not fix any collision;
therefore the fixed point index is $0$. 
Actually, what happens is that for some values of the coefficients
there are no fixed points, for some other values there are.

\subsection{Non-collinear solution}
Let $x_{1}=|z_2-z_3|$, $x_2 = |z_1-z_3|$ and $x_3=|z_1-z_2|$.
We have
\[
U = m_{12} \phi(|z_1 - z_2|) + m_{13} \phi( | z_1-z_3| ) +
m_{23} \phi(|z_2-z_3|),
\]
that is
\[
U = m_{12} \phi ( x_3) + m_{13} \phi (x_2) + m_{23}\phi(x_1),
\]
while
\begin{eqnarray*}
I  &=& m_1 |z_1|^2 + m_2 |z_2|^2 + m_3 |z_3|^2 = \\
&=&(m_1+m_2+m_3)^{-1}(m_1m_2 x_3^2 + m_1m_3 x_2^2 + m_2m_3 x_1^2).
\end{eqnarray*}
The function $u\from X \to \R$ therefore induces  function
$\bar u: X/C \to S^1$ such that $\bar u  \pi = u$,
where $\pi\from X \to \RP^2$ denotes the map
which sends $[z_1:z_2:z_3]\in X]$ to $[x_1:x_2:x_3]\in \RP^2$. 
Let $X_1$ denote the set of non-collinear configurations in $X$.
The image $\pi(X_1) \subset \RP^2$ is homeomorphic to 
the quotient $X_1/C$. Now, we can apply the same argument 
as in lemma \ref{lemma:2} and show that the critical points of $\bar u$
correspond to fixed points of the 
self-map $X_1/C \to \RP^2$ given by 
\[
\fbar [x_1:x_2:x_3] = [ \dfrac{m_{23}}{m_2m_3} x_1^{\alpha+1} :
 \dfrac{m_{13}}{m_1m_3} x_2^{\alpha+1} :
 \dfrac{m_{12}}{m_1m_2} x_3^{\alpha+1} ].
\]
Clearly, if any one of the $m_{ij}$ vanishes then 
there are no such fixed points in $X_1$ (a configuration
of type $[z:-z:0]$ is always collinear).
So in the affine chart the fixed points are the
compatible solutions
of the system of equations
\begin{eqnarray*}
x^{|\alpha|}&=&  \dfrac{m_{23} m_1}{m_{12} m_3} \\
y^{|\alpha|} &=& \dfrac{m_{13} m_2}{m_{12} m_3} \\
\end{eqnarray*}

If the $m_{ij}$ do not have the same sign, then there 
are no solutions.
In case of the Newtonian case we have $m_{ij}=m_im_j$, and hence 
the
solution is given as expected by 
the equilateral triangle ($x_1=x_2=x_3=1$).
Otherwise the solutions are either one  
(if the compatibility conditions on the sides
of a triangle are fulfilled: $x_i+x_j>x_k$ for 
each permutation $(i,j,k)$),
or none.

We can summarize the results in the following proposition.

\begin{propo}\label{propo:summa}
Given the coefficients $m_{12}$, $m_{13}$ and $m_{23}$,
masses $m_i>0$, $i=1,2,3$  and $\alpha<-1$,
if we define $\tilde m_1= m_1m_{23}$, 
$\tilde m_2 = m_2m_{13}$ and $\tilde m_3 = m_3m_{12}$,
the number of planar non-collinear central configurations  
$\#\CC_1/O(2)$ 
is:
\begin{enumerate}
\item 
$1$,  {if $\tilde m_{i}$ do not vanish, 
have the same sign, and
for every permutation $(i,j,k)$}
\[
(\tilde m_i)^{\frac{1}{|\alpha|}} +
(\tilde m_j)^{\frac{1}{|\alpha|}} >
(\tilde m_k)^{\frac{1}{|\alpha|}} 
\]
\item 
$0$, {otherwise}.
\end{enumerate}
The number of collinear central configurations is:
\begin{enumerate}
\item
If two of the $\tilde m_i$ vanish: $1$.
\item
If one of the $\tilde m_i$ vanishes and 
the sign of the left $\tilde m_{j}$ and $\tilde m_{k}$
is the same: at least $1$. 
\item
If $\tilde m_i \neq 0$ for $i=1,2,3$ and the $\tilde m_i$
have the same sign: $3$.
Otherwise, if the sign of the $\tilde m_i$ changes: at least $1$.
\end{enumerate}
\end{propo}

\section{Computing central configurations}
As a consequence of lemma \ref{lemma:2}, 
the central configurations are the solutions of 
the system of $2(n-2)$ equations in $2(n-2)$ unknown
variables. Let 
$A$ be the affine chart in $\CP^{n-2}$ of $z_{n-1}\neq 0$.
Let $x=(x_1,\dots, x_{n-2})$ be the affine coordinates.
In the chart $A$ the map $f$ is of the form
$x\mapsto w_i/w_{n-1}$,
where $w_i$ is the $i$-th component of $f$ and it is 
evaluated on the point
\[
[x_1:x_2:\dots : x_{n-2} : 1 : -1 -\sum_{i=1}^{n-2} {x_i}].
\]
Thus the $n-2$ equations are
\[
w_{n-1} x_i = w_i,\ \  (i=1,\dots, n-2). 
\]
I have computed the approximate solutions of such system 
for $n=3,\dots, 10$.
Unfortunately, only in the case of $n=4$ Albouy \cite{A95,AC98}
proved that these are indeed solutions and the only solutions
(actually, Morse equality \ref{eq:morse} below gives an alternate proof
of the fact that a non-collinear 
solution with symmetry different from $3$ and $4$ exists).
For $n\geq 5$ the fact that the system is solved up to the 
machine precision (in this case $10^{-15}$) does not imply
that the solutions actually exists. A computer-assisted 
proof was shown by Kotsireas \cite{Ko01}  for $n=4,5$.

A few words on the implementation: the program which 
computed the configurations in the appendix is written in FORTRAN 95,
using the SLATEC F77 library, and partially in the language of MAPLE
for post-processing. 
The algorithm is simple: 
a central configuration of equal masses can be ordered
in a way that $|z_1| \leq |z_2| \leq \dots \leq |z_n|$.
Therefore without loss of generality we can compute 
the configurations in the affine chart above. 
The root-finding subroutine starts by a random point
until reaches a solution or a failure (there are different reasons 
for the failure, like approaching a collision or too many steps
without progress).
Since the norms of $z_i$ are ordered, we can assume that 
and such 
that $|x_i|\leq 1$ for $i=1,\dots, n-2$.
Therefore the random starting point can be chosen in the cube
$[-1,1]^{2(n-2)}$.
After a solution is found, the program stores it in a list,
and computes the monic polynomial
\[
p(z) = \prod_{i=1}^n (z-z_i) = z^n + a_{n-2}z^{n-2} + \dots + a_1z  + a_0.
\]
Given two solutions, they coincide up to permuting the particles, 
scaling, rotation or 
conjugation if and only if  the coefficients $a_i$ and $a'_i$ 
of the corresponding 
polynomials are related by the equations
\[
b^{i}a_i = a'_i
\]
or 
\[
b^ia_i = \bar a'_i,
\]
where $b$ is a suitable complex number.
Then, if the configuration is not in the list of known configurations,
it is added to it. After the program cannot find any new 
configuration for a long enough time (depending on the number of masses),
the list of configurations found is piped to a MAPLE filter that 
cross-checks the solutions, computes the 
reduced potential $U\sqrt{I}$, the
critical point  index,
the fixed point index and the order of the isotropy group of the configuration
(the symmetric group $\Sigma_n$ acts on $X$, therefore each solution 
has a isotropy group with respect to this action: in the tables
we count this isotropy, in general; for a few configurations
which are not symmetric with respect to a reflection, we 
denote the isotropy $1/2$).

The time spent by the program in computing the solutions is 
a tiny part, whenever compared to the time spent in searching for 
solutions once the list is complete. Unfortunately, not only these 
numerical experiments do not imply that the approximate solutions 
are solutions, but there is no proof of the fact that there are no
extra solutions. Some hope that the list, at least for small $n$,
is complete, comes from the Morse equality
\begin{equation}\label{eq:morse}
\sum_{k} (-1)^k\nu_k = \chi(X) = (-1)^n(n-2)!.
\end{equation}
Let $i(z)$ denote the size of the isotropy group of 
the central configuration $z\in \CC_1/SO(2)$ with respect 
to the action of $\Sigma_n$, and 
$h(z)$ its Morse index.
A consequence of equation \ref{eq:morse} 
is 
the following equation.
\begin{equation}\label{eq:my}
\sum_{z\in \CC_1/SO(2)} \dfrac{(-1)^{h(z)}}{i(z)}  = \dfrac{(-1)^n}{n(n-1)}.
\end{equation}
Therefore we can test equation \ref{eq:my} against the collected 
data, keeping in mind that in the tables the 
configurations are listed in $\CC_1/O(2)$, hence if there is 
no axis of symmetry the configuration contributes to the sum twice.
For example, for $n=4$ we have 
\[
\frac{1}{4} - 1 + \frac{1}{3} + \frac{1}{2} = \frac{1}{12};
\]
for $n=5$,
\[
\frac{1}{5} + \frac{1}{4} - 1 + 1 - \frac{1}{2} = -\frac{1}{20};
\]
for $n=6$,
\[
\frac{1}{5} + \frac{1}{3} -\frac{1}{6} 
-1
+\frac{1}{2}
+1
-\frac{1}{3} -1 + \frac{1}{2} = \frac{1}{30};
\]
for $n=7$ analogously we obtain $-\frac{1}{42}$;
for $n=8$ we have the first examples of configurations without 
an axis of symmetry: they contribute with $\dfrac{1}{1/2} = 2$ 
instead of $1$, and the sum gives $\frac{1}{56}$.
For $n=9$, it gives the expected $-\dfrac{1}{72}$.
Unfortunately, for $n=10$ the sum is 
$-\dfrac{269}{90}$, thus there is a missing term $+3$.
It is likely that the computation of the isotropy was broken
somewhere or simply that there are some configurations missing.

\section{Remarks}

\begin{remark}
In case of the logarithmic potential, the collinear central configurations
are the zeroes of the Hermite polynomials (\cite{Sz75},
Theorem 6.7.3 pag.~141). This result can be in some
sense extended
to the planar central configurations as follows: the planar central
configurations of $n$ bodies
are the zeroes of degree $n$ polynomials $p(x)$ with complex
coefficients such that $p(x) = 0$ $\implies$
$ p''(x) + \lambda \overline{x} p'(x)=0$,
with $\lambda$ real non-zero parameter.
That is, $p(x)$ is a polynomial of degree $n$
satisfying  the differential equation
$p'' + \lambda \overline{x} p' + c p =0$,
where $c$ is a suitable function of $x$.
\end{remark}

\begin{remark}
It is apparent from the data shown below 
that the fixed point index of 
a central configuration is likely to be equal to minus one to the power
the critical point index. This is certainly true for 
the gradient map, but $f$ is not equal (nor homotopic) to the gradient map.
This relation is true for collinear configurations (embedded in the plane),
and it is likely that it is true for all the central configurations.
Unfortunately I do not know any proof of it.
\end{remark}

\begin{remark}
We see from the examples in the tables that if a central configuration 
has a rotational symmetry, then it has a symmetry axis.
Thus the symmetry group of the configuration is 
either trivial or a dihedral group.
In the tables there are no central configurations with a pure rotation symmetry.
Is it true in general?
\end{remark}


\part*{Appendix}
\section*{Tables of central configurations}
{\tiny
\input{nbody_3}
\input{nbody_4}
\input{nbody_5}
\input{nbody_6}
\input{nbody_7}
\input{nbody_8}
\input{nbody_9}
\input{nbody_10}
}

\end{document}

%% file: nbody_3.tex
\begin{longtable}{|ll|}
\hline
\multicolumn{2}{|c|}{\emph{Central configurations for 3 bodies}}\\ 
\hline
\hline \begin{minipage}{0.4\textwidth}
\begin{center} 

\includegraphics[width=\textwidth]{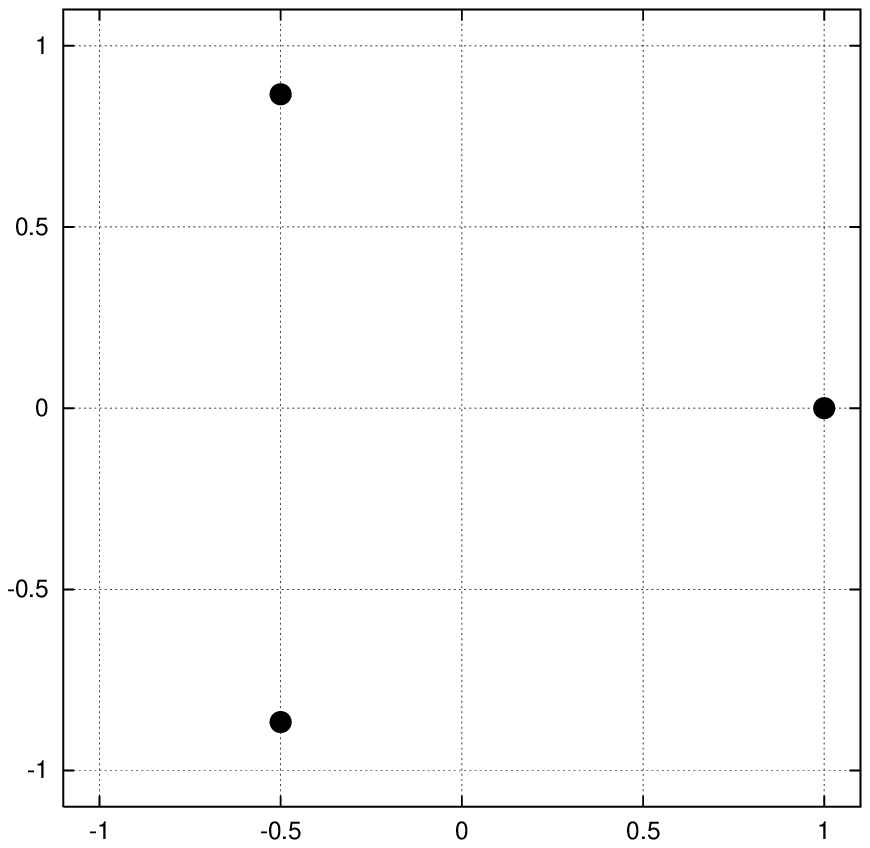} 
\end{center}\end{minipage}
&
\fbox{\begin{minipage}[t]{0.6\textwidth}
\begin{tabular}{r@{ $ = $ }l}
$U\sqrt{I}$ &    $3.00000000 $ \\  
Crit ind & $0$ \\ 
FP ind & $1$ \\ 
Isotropy & $3$\\ 
$z_1$ & $(-.50000000,.86602540)$\\
$z_2$ & $(1.00000000,0.00000000)$\\
$z_3$ & $(-.50000000,-.86602540)$\\
\end{tabular}\end{minipage}}
\\ 
\hline \begin{minipage}{0.4\textwidth}
\begin{center} 

\includegraphics[width=\textwidth]{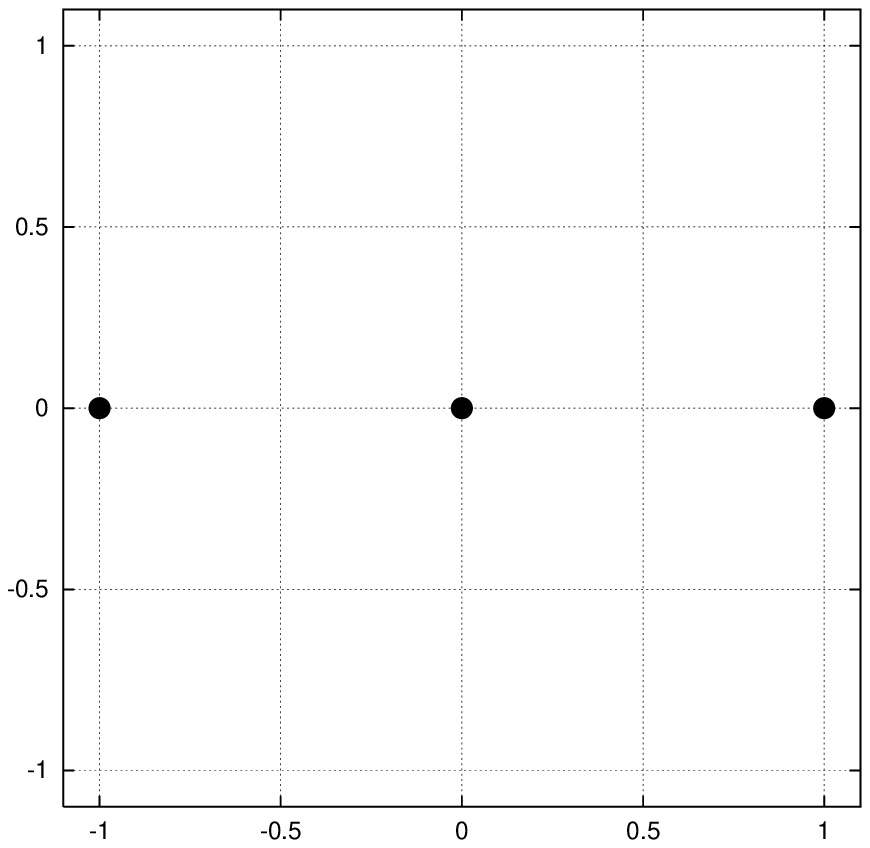} 
\end{center}\end{minipage}
&
\fbox{\begin{minipage}[t]{0.6\textwidth}
\begin{tabular}{r@{ $ = $ }l}
$U\sqrt{I}$ &    $3.53553391 $ \\  
Crit ind & $1$ \\ 
FP ind & $-1$ \\ 
Isotropy & $2$\\ 
$z_1$ & $(0.00000000,0.00000000)$\\
$z_2$ & $(1.00000000,0.00000000)$\\
$z_3$ & $(-1.00000000,0.00000000)$\\
\end{tabular}\end{minipage}}
\\ 
\hline\hline
\end{longtable}

%% file: nbody_4.tex
\begin{longtable}{|ll|}
\hline
\multicolumn{2}{|c|}{\emph{Central configurations for 4 bodies}}\\ 
\hline
\hline \begin{minipage}{0.4\textwidth}
\begin{center} 

\includegraphics[width=\textwidth]{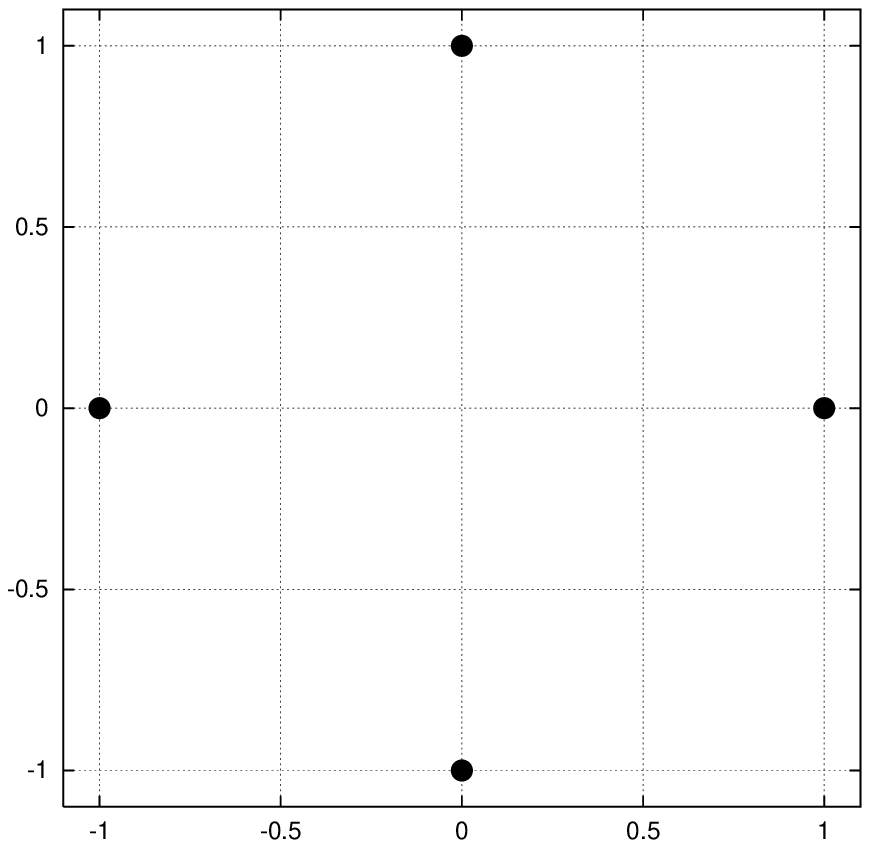} 
\end{center}\end{minipage}
&
\fbox{\begin{minipage}[t]{0.6\textwidth}
\begin{tabular}{r@{ $ = $ }l}
$U\sqrt{I}$ &    $7.65685425 $ \\  
Crit ind & $0$ \\ 
FP ind & $1$ \\ 
Isotropy & $4$\\ 
$z_1$ & $(.00000000,1.00000000)$\\
$z_2$ & $(.00000000,-1.00000000)$\\
$z_3$ & $(1.00000000,0.00000000)$\\
$z_4$ & $(-1.00000000,0.00000000)$\\
\end{tabular}\end{minipage}}
\\ 
\hline \begin{minipage}{0.4\textwidth}
\begin{center} 

\includegraphics[width=\textwidth]{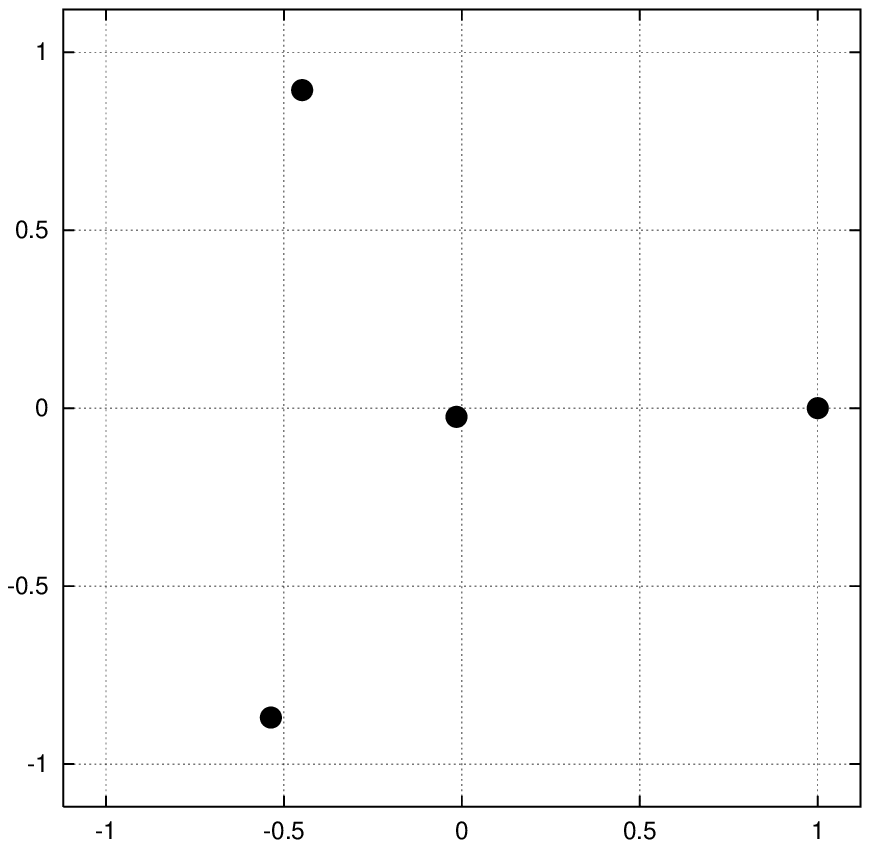} 
\end{center}\end{minipage}
&
\fbox{\begin{minipage}[t]{0.6\textwidth}
\begin{tabular}{r@{ $ = $ }l}
$U\sqrt{I}$ &    $8.19608063 $ \\  
Crit ind & $1$ \\ 
FP ind & $-1$ \\ 
Isotropy & $1$\\ 
$z_1$ & $(-.01495777,-.02424649)$\\
$z_2$ & $(-.44867582,.89369458)$\\
$z_3$ & $(1.00000000,0.00000000)$\\
$z_4$ & $(-.53636641,-.86944809)$\\
\end{tabular}\end{minipage}}
\\ 
\hline \begin{minipage}{0.4\textwidth}
\begin{center} 

\includegraphics[width=\textwidth]{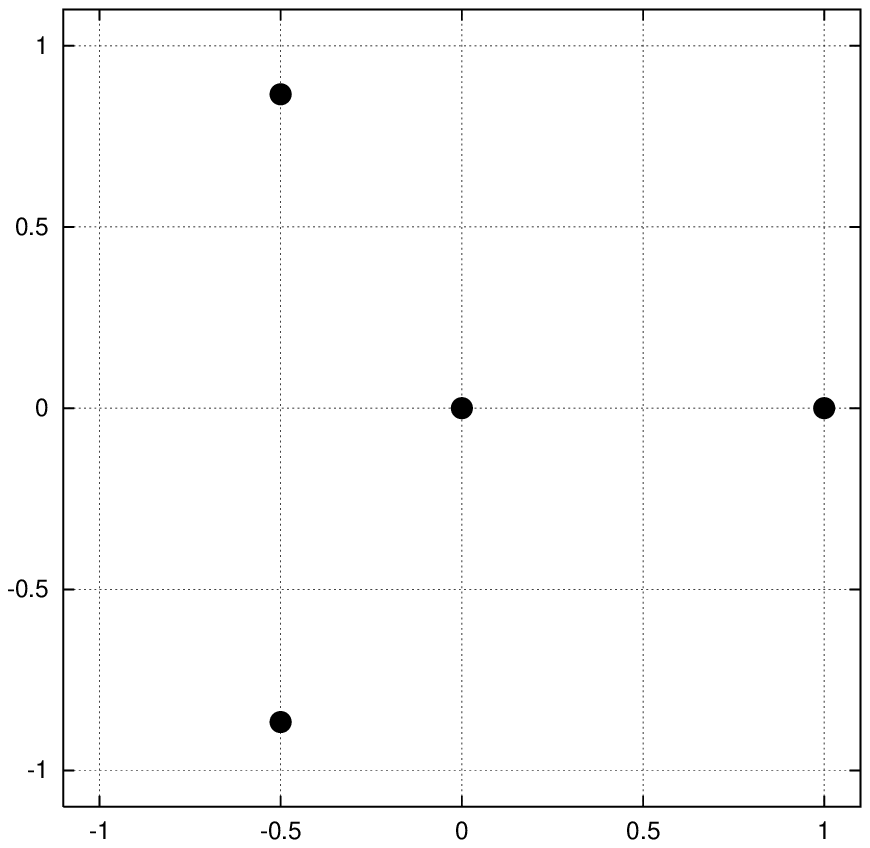} 
\end{center}\end{minipage}
&
\fbox{\begin{minipage}[t]{0.6\textwidth}
\begin{tabular}{r@{ $ = $ }l}
$U\sqrt{I}$ &    $8.19615242 $ \\  
Crit ind & $2$ \\ 
FP ind & $1$ \\ 
Isotropy & $3$\\ 
$z_1$ & $(-.00000000,.00000000)$\\
$z_2$ & $(-.50000000,-.86602540)$\\
$z_3$ & $(1.00000000,0.00000000)$\\
$z_4$ & $(-.50000000,.86602540)$\\
\end{tabular}\end{minipage}}
\\
\hline\hline
\end{longtable}

%% file: nbody_5.tex
\begin{longtable}{|ll|}
\hline
\multicolumn{2}{|c|}{\emph{Some central configurations for 5 bodies}}\\ 
\hline
\hline \begin{minipage}{0.4\textwidth}
\begin{center} 

\includegraphics[width=\textwidth]{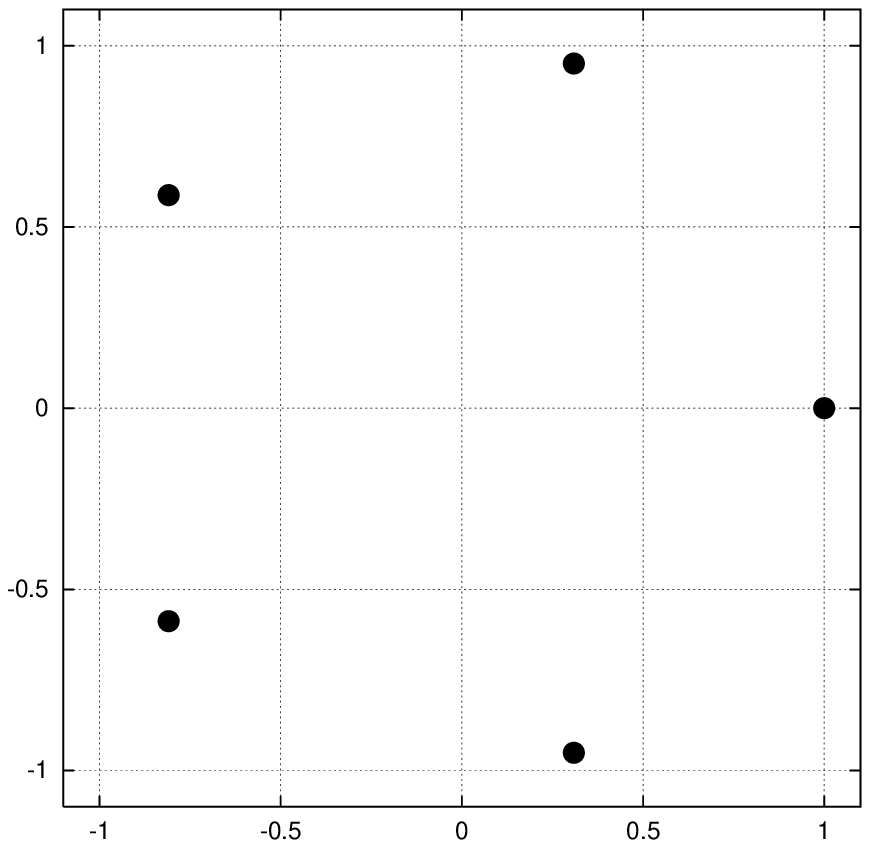} 
\end{center}\end{minipage}
&
\fbox{\begin{minipage}[t]{0.6\textwidth}
\begin{tabular}{r@{ $ = $ }l}
$U\sqrt{I}$ &    $15.38841769 $ \\  
Crit ind & $0$ \\ 
FP ind & $1$ \\ 
Isotropy & $5$\\ 
$z_1$ & $(-.80901699,-.58778525)$\\
$z_2$ & $(.30901699,.95105652)$\\
$z_3$ & $(.30901699,-.95105652)$\\
$z_4$ & $(1.00000000,0.00000000)$\\
$z_5$ & $(-.80901699,.58778525)$\\
\end{tabular}\end{minipage}}
\\ 
\hline \begin{minipage}{0.4\textwidth}
\begin{center} 

\includegraphics[width=\textwidth]{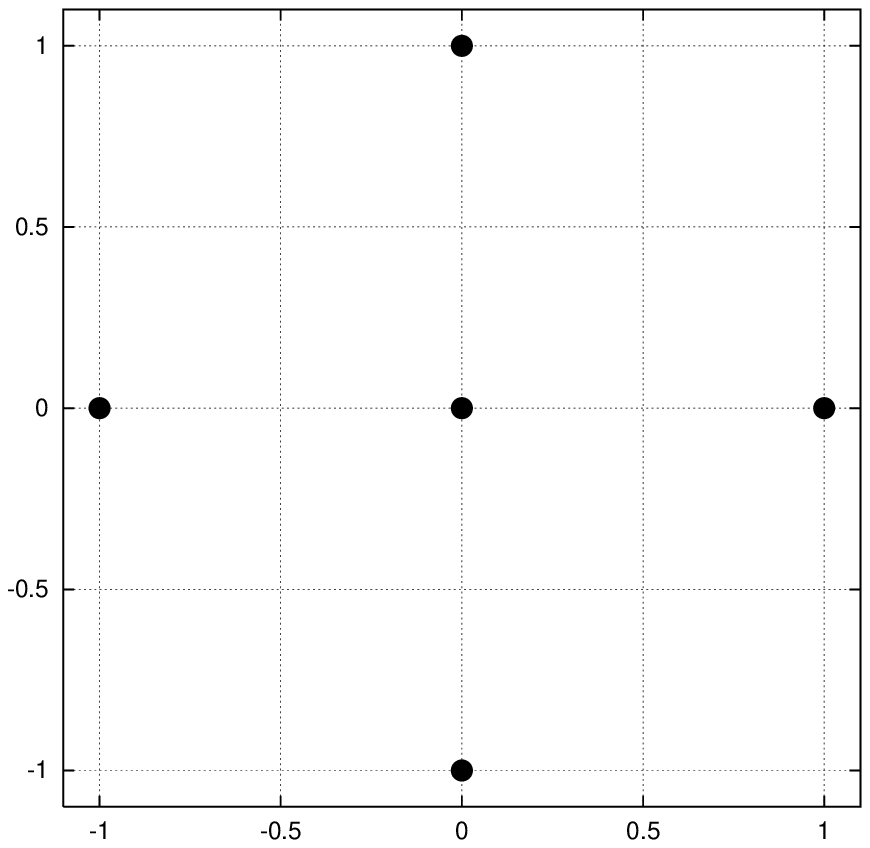} 
\end{center}\end{minipage}
&
\fbox{\begin{minipage}[t]{0.6\textwidth}
\begin{tabular}{r@{ $ = $ }l}
$U\sqrt{I}$ &    $15.65685425 $ \\  
Crit ind & $0$ \\ 
FP ind & $1$ \\ 
Isotropy & $4$\\ 
$z_1$ & $(.00000000,-.00000000)$\\
$z_2$ & $(-1.00000000,.00000000)$\\
$z_3$ & $(-.00000000,-1.00000000)$\\
$z_4$ & $(1.00000000,0.00000000)$\\
$z_5$ & $(.00000000,1.00000000)$\\
\end{tabular}\end{minipage}}
\\ 
\hline \begin{minipage}{0.4\textwidth}
\begin{center} 

\includegraphics[width=\textwidth]{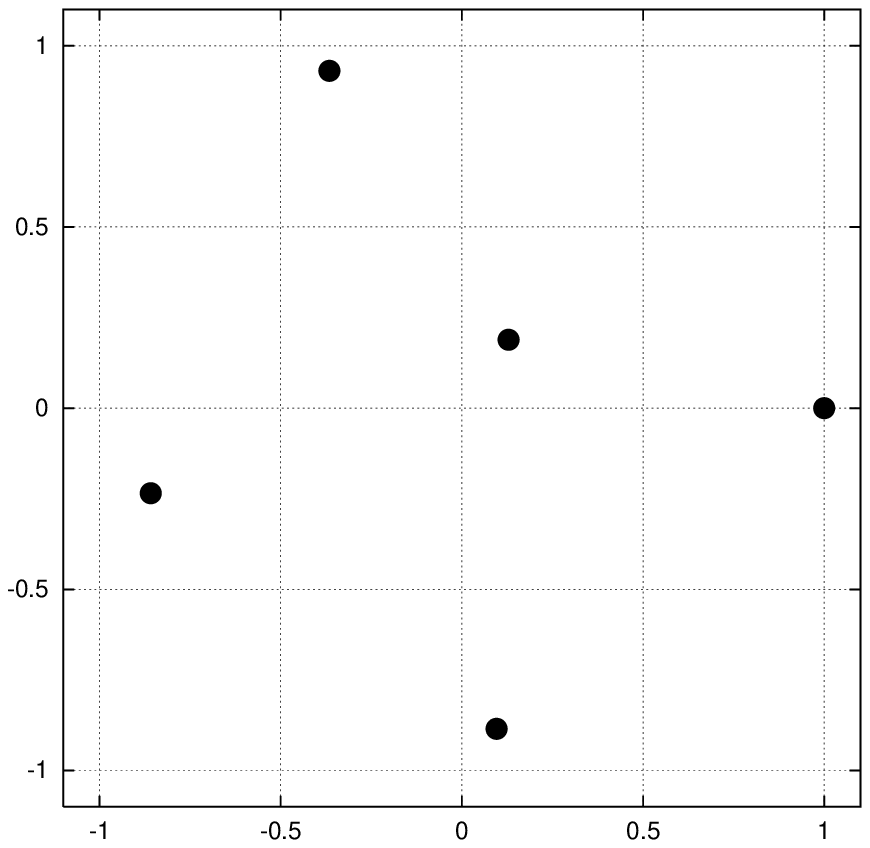} 
\end{center}\end{minipage}
&
\fbox{\begin{minipage}[t]{0.6\textwidth}
\begin{tabular}{r@{ $ = $ }l}
$U\sqrt{I}$ &    $15.68397123 $ \\  
Crit ind & $1$ \\ 
FP ind & $-1$ \\ 
Isotropy & $1$\\ 
$z_1$ & $(.12867234,.18881180)$\\
$z_2$ & $(.09578559,-.88501108)$\\
$z_3$ & $(-.85873077,-.23452283)$\\
$z_4$ & $(1.00000000,0.00000000)$\\
$z_5$ & $(-.36572717,.93072211)$\\
\end{tabular}\end{minipage}}
\\
\hline\hline
\end{longtable}

%% file: nbody_6.tex
\begin{longtable}{|ll|}
\hline
\multicolumn{2}{|c|}{\emph{Some central configurations for 6 bodies}}\\ 
\hline
\hline \begin{minipage}{0.4\textwidth}
\begin{center} 

\includegraphics[width=\textwidth]{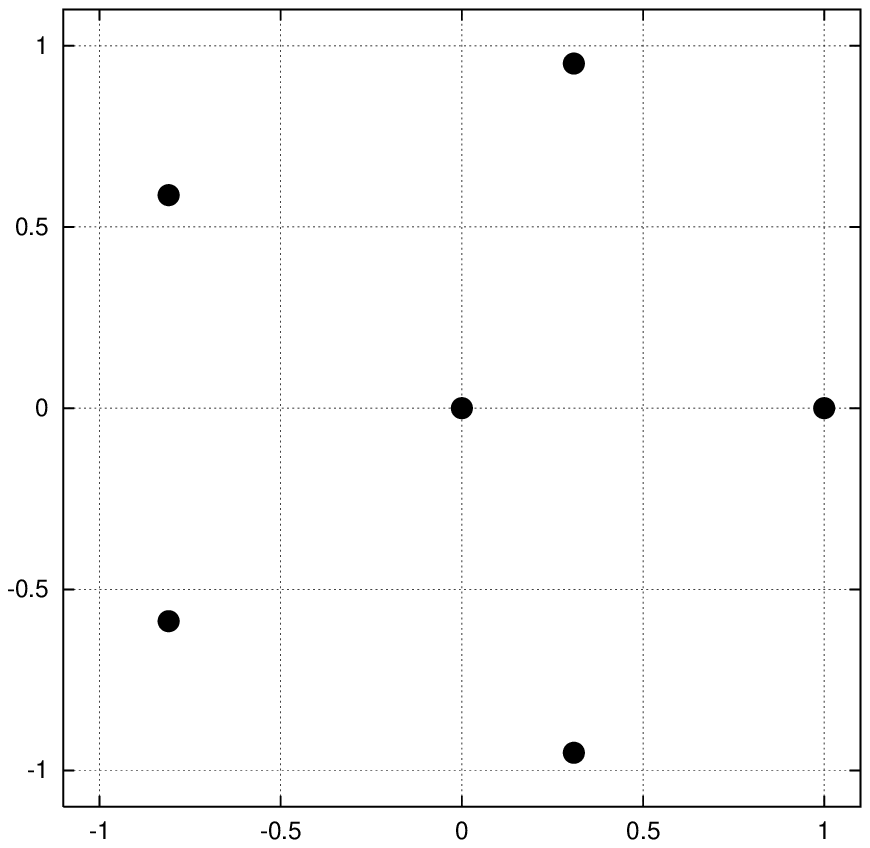} 
\end{center}\end{minipage}
&
\fbox{\begin{minipage}[t]{0.6\textwidth}
\begin{tabular}{r@{ $ = $ }l}
$U\sqrt{I}$ &    $26.56875757 $ \\  
Crit ind & $0$ \\ 
FP ind & $1$ \\ 
Isotropy & $5$\\ 
$z_1$ & $(-.00000000,.00000000)$\\
$z_2$ & $(.30901699,-.95105652)$\\
$z_3$ & $(-.80901699,.58778525)$\\
$z_4$ & $(-.80901699,-.58778525)$\\
$z_5$ & $(1.00000000,0.00000000)$\\
$z_6$ & $(.30901699,.95105652)$\\
\end{tabular}\end{minipage}}
\\ 
\hline \begin{minipage}{0.4\textwidth}
\begin{center} 

\includegraphics[width=\textwidth]{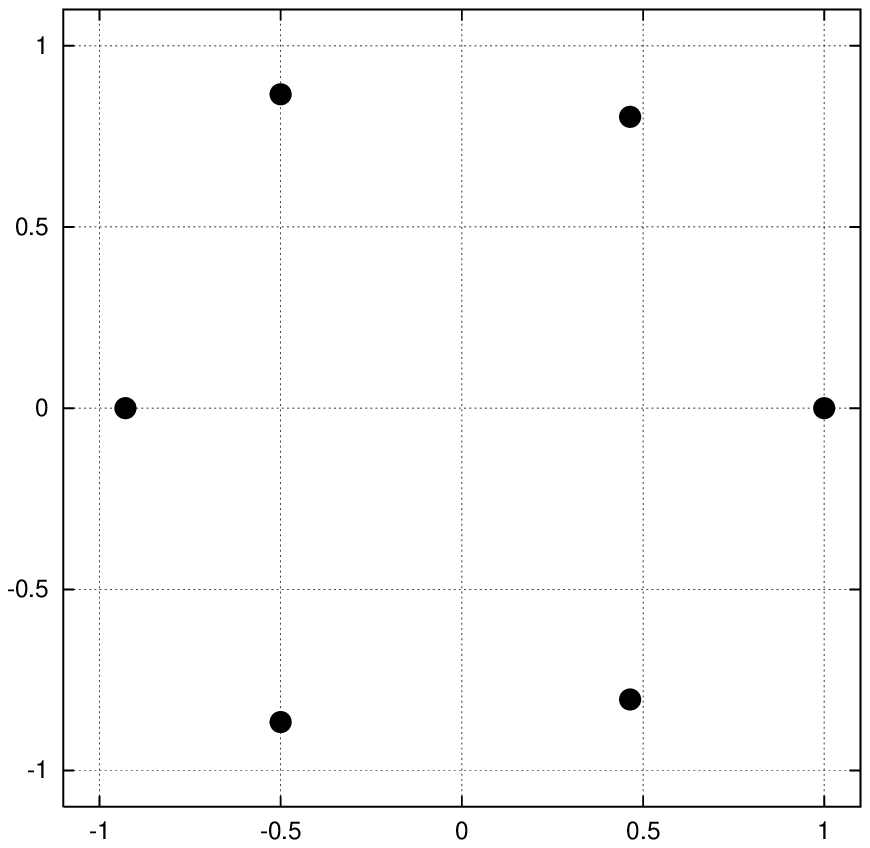} 
\end{center}\end{minipage}
&
\fbox{\begin{minipage}[t]{0.6\textwidth}
\begin{tabular}{r@{ $ = $ }l}
$U\sqrt{I}$ &    $26.85636352 $ \\  
Crit ind & $0$ \\ 
FP ind & $1$ \\ 
Isotropy & $3$\\ 
$z_1$ & $(.46417852,.80398078)$\\
$z_2$ & $(-.92835704,-.00000000)$\\
$z_3$ & $(.46417852,-.80398078)$\\
$z_4$ & $(-.50000000,-.86602540)$\\
$z_5$ & $(1.00000000,0.00000000)$\\
$z_6$ & $(-.50000000,.86602540)$\\
\end{tabular}\end{minipage}}
\\ 
\hline \begin{minipage}{0.4\textwidth}
\begin{center} 

\includegraphics[width=\textwidth]{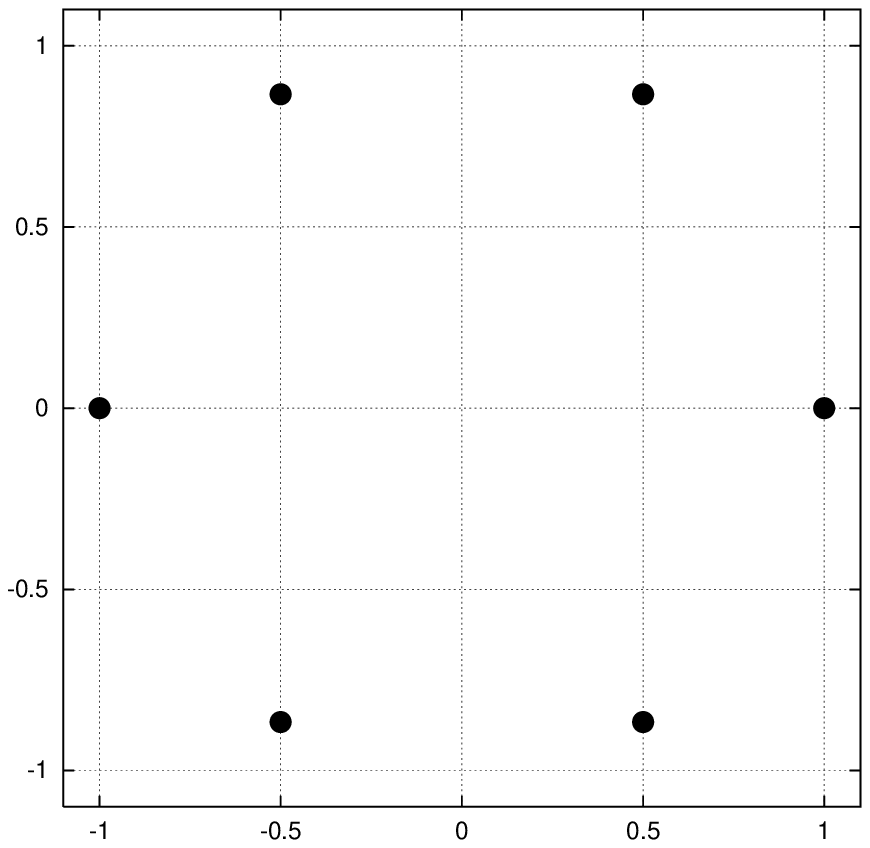} 
\end{center}\end{minipage}
&
\fbox{\begin{minipage}[t]{0.6\textwidth}
\begin{tabular}{r@{ $ = $ }l}
$U\sqrt{I}$ &    $26.85645445 $ \\  
Crit ind & $1$ \\ 
FP ind & $-1$ \\ 
Isotropy & $6$\\ 
$z_1$ & $(.50000000,-.86602540)$\\
$z_2$ & $(-1.00000000,-.00000000)$\\
$z_3$ & $(.50000000,.86602540)$\\
$z_4$ & $(-.50000000,.86602540)$\\
$z_5$ & $(1.00000000,0.00000000)$\\
$z_6$ & $(-.50000000,-.86602540)$\\
\end{tabular}\end{minipage}}
\\
\hline\hline
\end{longtable}

%% file: nbody_7.tex
\begin{longtable}{|ll|}
\hline
\multicolumn{2}{|c|}{\emph{Some central configurations for 7 bodies}}\\ 
\hline
\hline \begin{minipage}{0.4\textwidth}
\begin{center} 

\includegraphics[width=\textwidth]{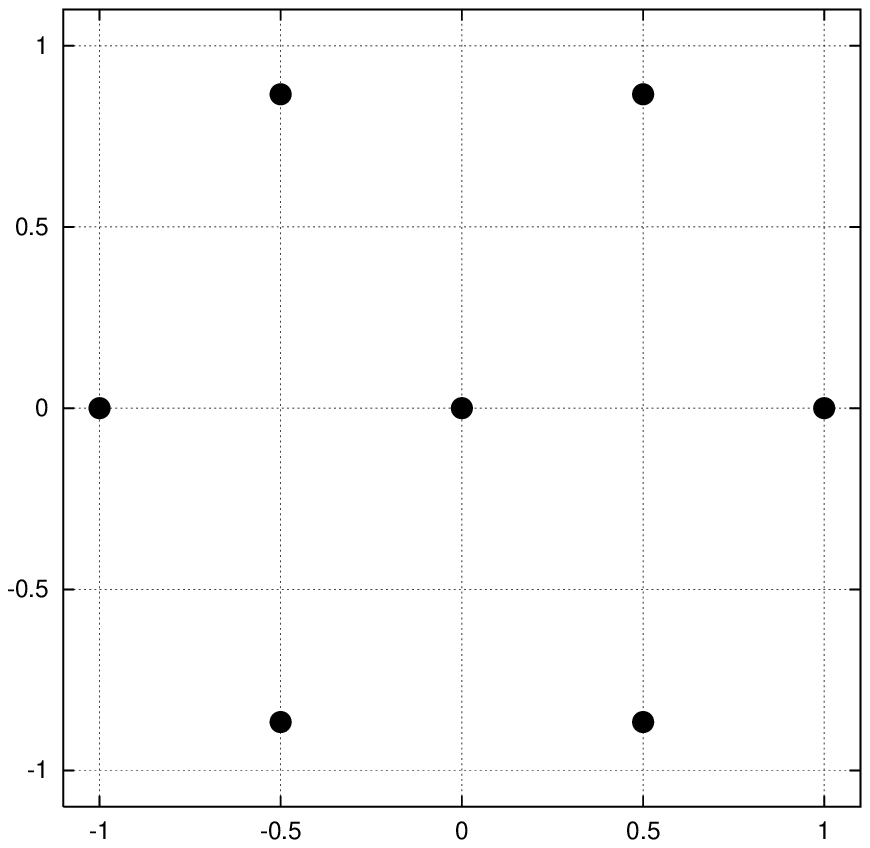} 
\end{center}\end{minipage}
&
\fbox{\begin{minipage}[t]{0.6\textwidth}
\begin{tabular}{r@{ $ = $ }l}
$U\sqrt{I}$ &    $41.55339290 $ \\  
Crit ind & $0$ \\ 
FP ind & $1$ \\ 
Isotropy & $6$\\ 
$z_1$ & $(-.00000000,-.00000000)$\\
$z_2$ & $(-.50000000,-.86602540)$\\
$z_3$ & $(.50000000,-.86602540)$\\
$z_4$ & $(.50000000,.86602540)$\\
$z_5$ & $(-1.00000000,.00000000)$\\
$z_6$ & $(1.00000000,0.00000000)$\\
$z_7$ & $(-.50000000,.86602540)$\\
\end{tabular}\end{minipage}}
\\ 
\hline \begin{minipage}{0.4\textwidth}
\begin{center} 

\includegraphics[width=\textwidth]{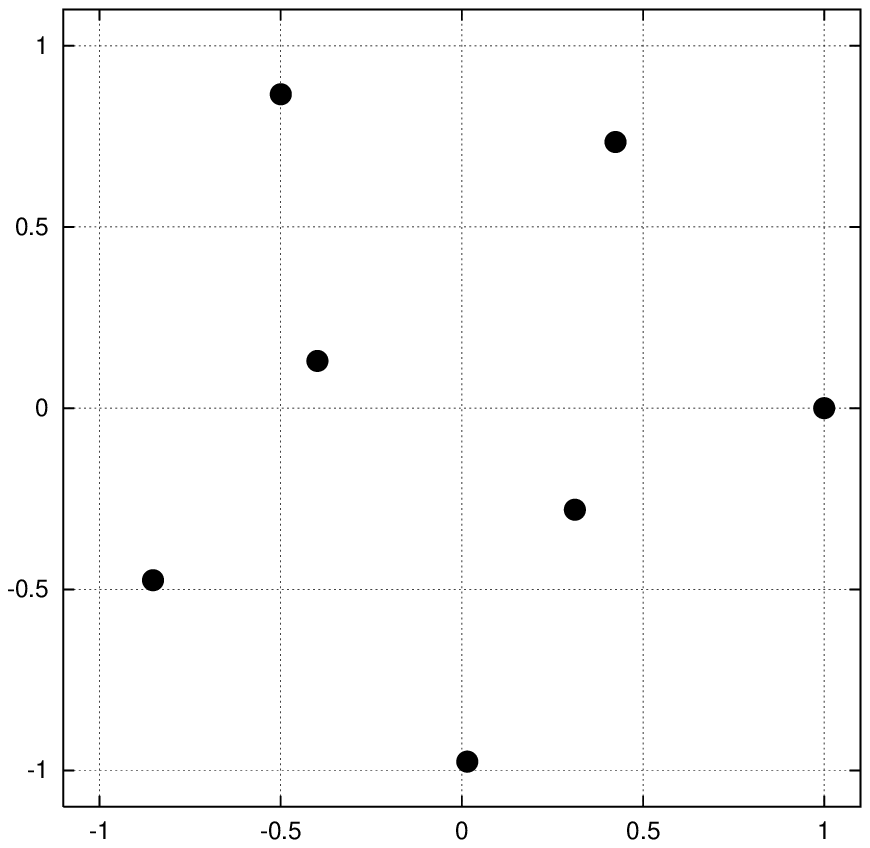} 
\end{center}\end{minipage}
&
\fbox{\begin{minipage}[t]{0.6\textwidth}
\begin{tabular}{r@{ $ = $ }l}
$U\sqrt{I}$ &    $42.32859664 $ \\  
Crit ind & $1$ \\ 
FP ind & $-1$ \\ 
Isotropy & $1$\\ 
$z_1$ & $(.31194658,-.28006242)$\\
$z_2$ & $(-.39849482,.13018261)$\\
$z_3$ & $(.42399335,.73425006)$\\
$z_4$ & $(-.85245018,-.47481689)$\\
$z_5$ & $(.01487433,-.97565421)$\\
$z_6$ & $(1.00000000,0.00000000)$\\
$z_7$ & $(-.49986928,.86610086)$\\
\end{tabular}\end{minipage}}
\\ 
\hline \begin{minipage}{0.4\textwidth}
\begin{center} 

\includegraphics[width=\textwidth]{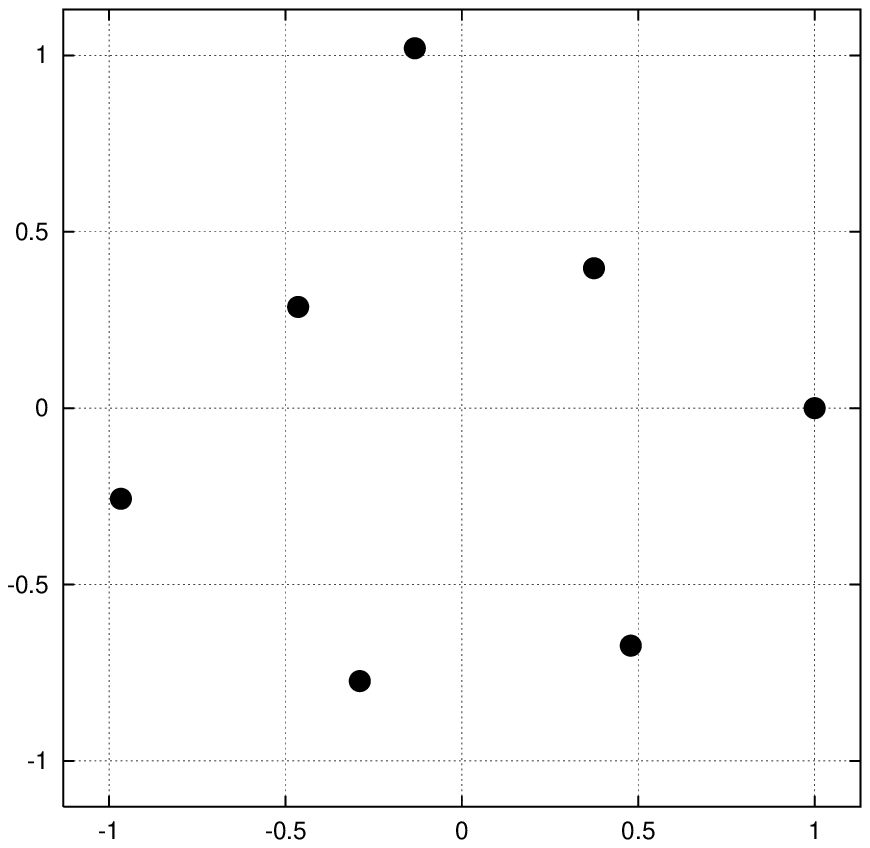} 
\end{center}\end{minipage}
&
\fbox{\begin{minipage}[t]{0.6\textwidth}
\begin{tabular}{r@{ $ = $ }l}
$U\sqrt{I}$ &    $42.45998988 $ \\  
Crit ind & $1$ \\ 
FP ind & $-1$ \\ 
Isotropy & $1$\\ 
$z_1$ & $(-.46413573,.28708218)$\\
$z_2$ & $(.37471387,.39677156)$\\
$z_3$ & $(.47873707,-.67327333)$\\
$z_4$ & $(-.28952417,-.77373246)$\\
$z_5$ & $(-.96637766,-.25712684)$\\
$z_6$ & $(1.00000000,0.00000000)$\\
$z_7$ & $(-.13341338,1.02027888)$\\
\end{tabular}\end{minipage}}
\\
\hline\hline
\end{longtable}

%% file: nbody_8.tex
\begin{longtable}{|ll|}
\hline
\multicolumn{2}{|c|}{\emph{Some central configurations for 8 bodies}}\\ 
\hline
\hline \begin{minipage}{0.4\textwidth}
\begin{center} 

\includegraphics[width=\textwidth]{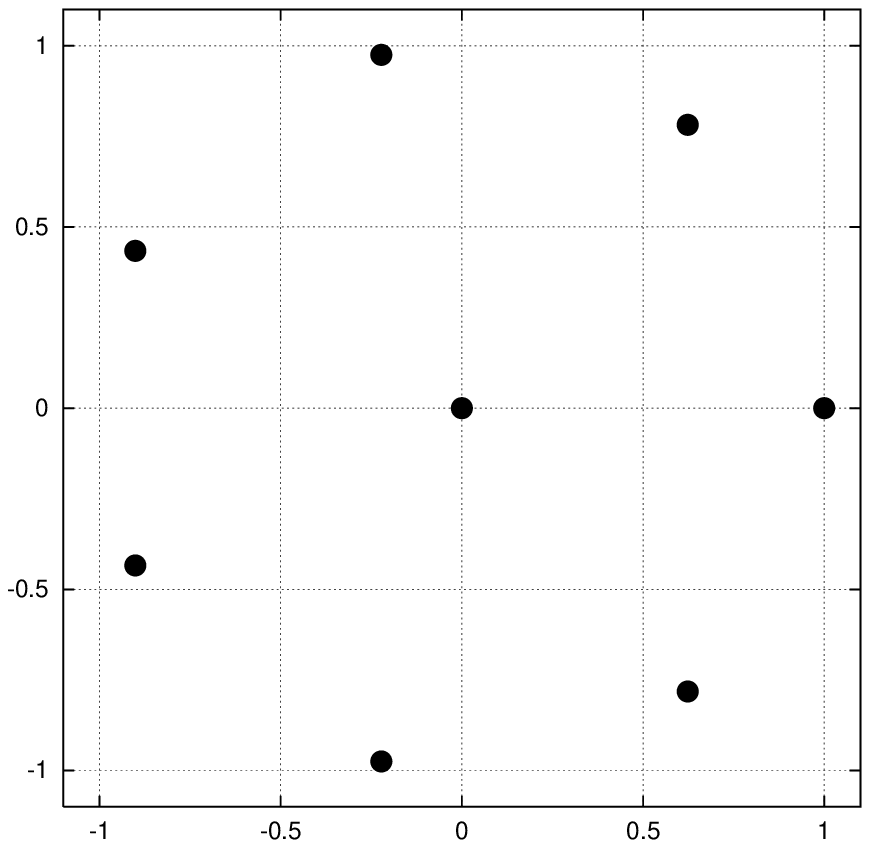} 
\end{center}\end{minipage}
&
\fbox{\begin{minipage}[t]{0.6\textwidth}
\begin{tabular}{r@{ $ = $ }l}
$U\sqrt{I}$ &    $61.20510193 $ \\  
Crit ind & $0$ \\ 
FP ind & $1$ \\ 
Isotropy & $7$\\ 
$z_1$ & $(.00000000,.00000000)$\\
$z_2$ & $(.62348980,.78183148)$\\
$z_3$ & $(-.90096887,.43388374)$\\
$z_4$ & $(-.22252093,-.97492791)$\\
$z_5$ & $(-.90096887,-.43388374)$\\
$z_6$ & $(.62348980,-.78183148)$\\
$z_7$ & $(1.00000000,0.00000000)$\\
$z_8$ & $(-.22252093,.97492791)$\\
\end{tabular}\end{minipage}}
\\ 
\hline \begin{minipage}{0.4\textwidth}
\begin{center} 

\includegraphics[width=\textwidth]{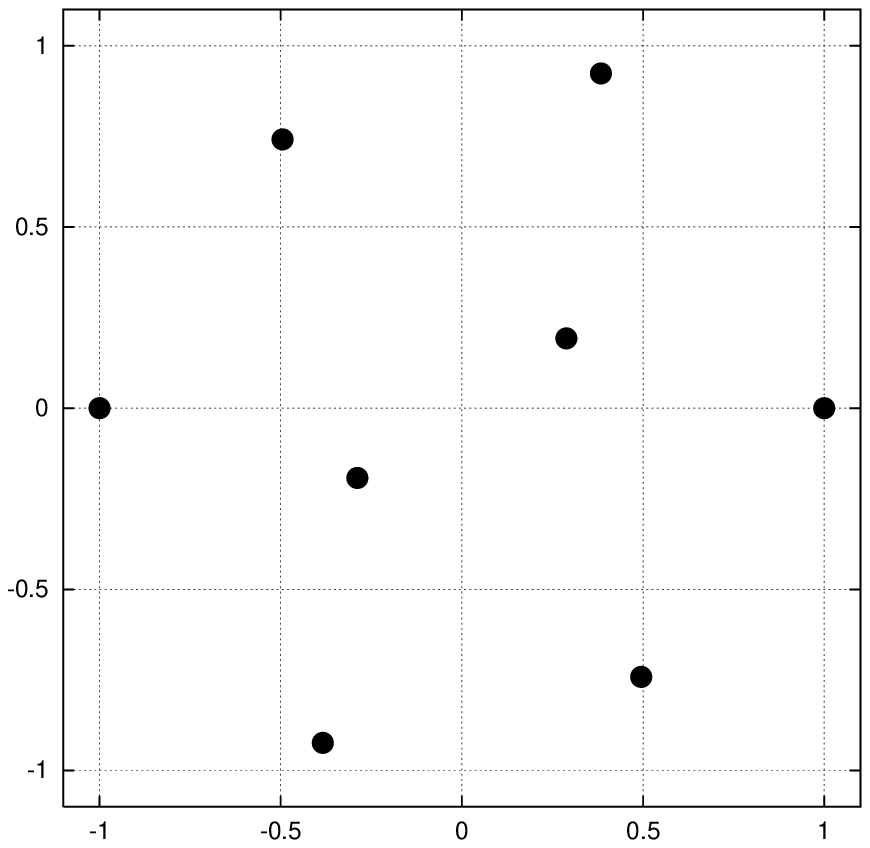} 
\end{center}\end{minipage}
&
\fbox{\begin{minipage}[t]{0.6\textwidth}
\begin{tabular}{r@{ $ = $ }l}
$U\sqrt{I}$ &    $61.66790403 $ \\  
Crit ind & $1$ \\ 
FP ind & $-1$ \\ 
Isotropy & $2$\\ 
$z_1$ & $(-.28827131,-.19240111)$\\
$z_2$ & $(.28827131,.19240111)$\\
$z_3$ & $(-.49497029,.74160556)$\\
$z_4$ & $(.49497029,-.74160556)$\\
$z_5$ & $(-.38363904,-.92348313)$\\
$z_6$ & $(-1.00000000,-.00000000)$\\
$z_7$ & $(1.00000000,0.00000000)$\\
$z_8$ & $(.38363904,.92348313)$\\
\end{tabular}\end{minipage}}
\\
\hline\hline
\end{longtable}

%% file: nbody_9.tex
\begin{longtable}{|ll|}
\hline
\multicolumn{2}{|c|}{\emph{Some central configurations for 9 bodies}}\\ 
\hline
\hline \begin{minipage}{0.4\textwidth}
\begin{center} 

\includegraphics[width=\textwidth]{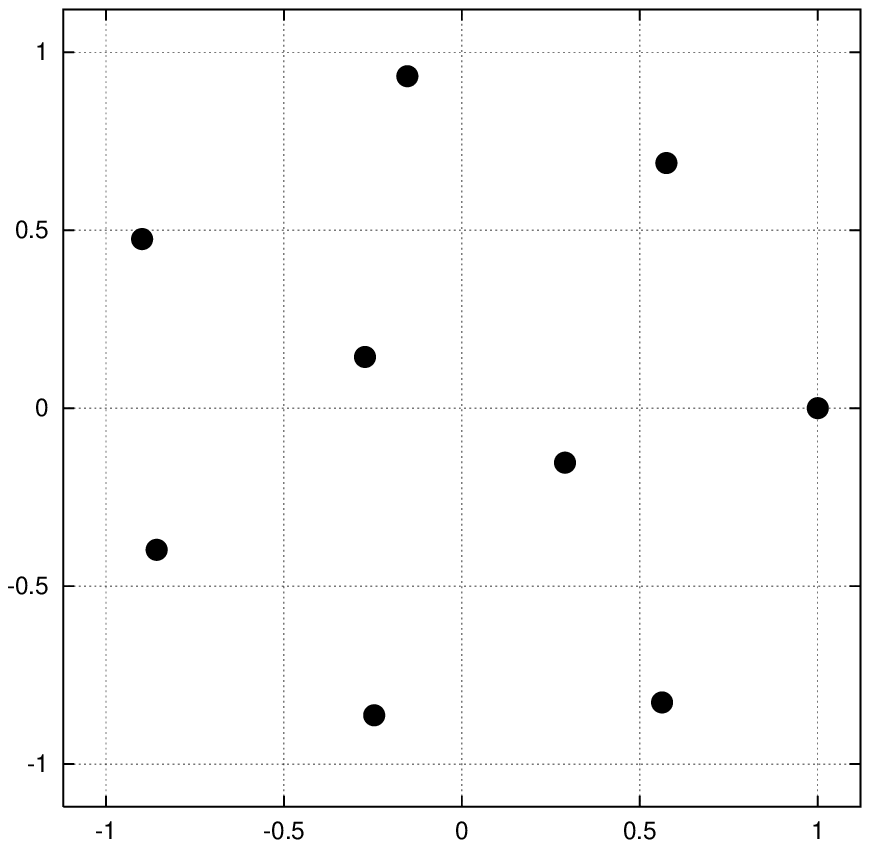} 
\end{center}\end{minipage}
&
\fbox{\begin{minipage}[t]{0.6\textwidth}
\begin{tabular}{r@{ $ = $ }l}
$U\sqrt{I}$ &    $85.90069864 $ \\  
Crit ind & $0$ \\ 
FP ind & $1$ \\ 
Isotropy & $1$\\ 
$z_1$ & $(-.27202755,.14390992)$\\
$z_2$ & $(.28973718,-.15327879)$\\
$z_3$ & $(-.24618043,-.86277978)$\\
$z_4$ & $(.57473480,.68896565)$\\
$z_5$ & $(-.15320370,.93279998)$\\
$z_6$ & $(-.85733689,-.39819672)$\\
$z_7$ & $(.56265924,-.82668893)$\\
$z_8$ & $(1.00000000,0.00000000)$\\
$z_9$ & $(-.89838264,.47526867)$\\
\end{tabular}\end{minipage}}
\\ 
\hline \begin{minipage}{0.4\textwidth}
\begin{center} 

\includegraphics[width=\textwidth]{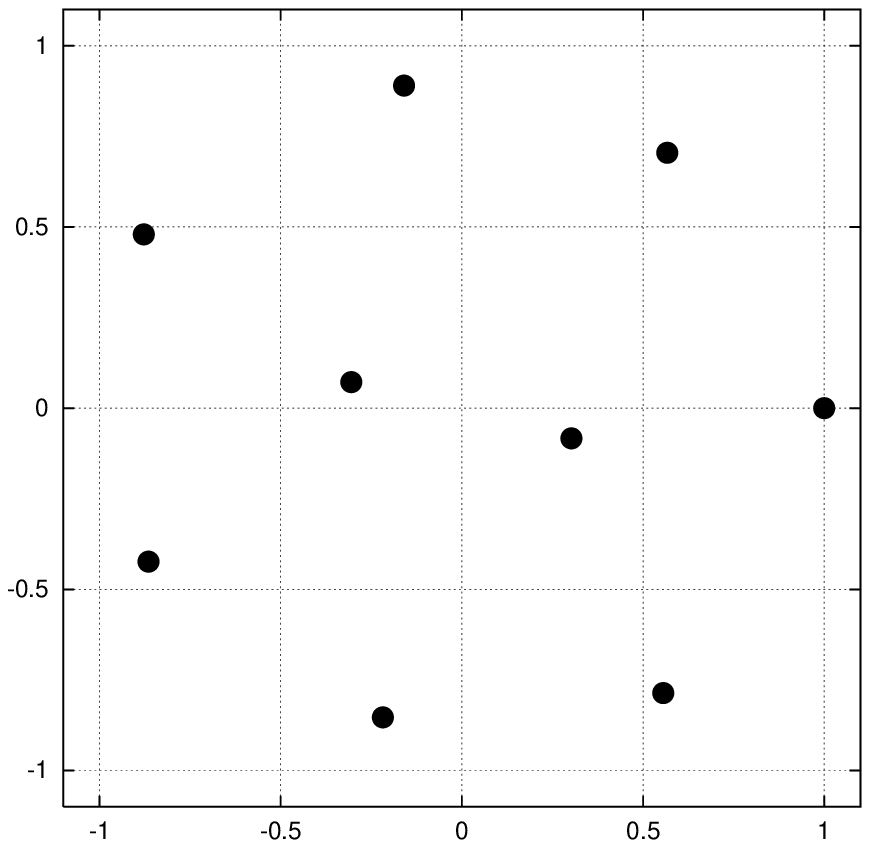} 
\end{center}\end{minipage}
&
\fbox{\begin{minipage}[t]{0.6\textwidth}
\begin{tabular}{r@{ $ = $ }l}
$U\sqrt{I}$ &    $85.90099418 $ \\  
Crit ind & $1$ \\ 
FP ind & $-1$ \\ 
Isotropy & $1$\\ 
$z_1$ & $(-.30525333,.07170100)$\\
$z_2$ & $(.30225296,-.08344927)$\\
$z_3$ & $(-.21785509,-.85303320)$\\
$z_4$ & $(.56690970,.70454393)$\\
$z_5$ & $(-.15965544,.89010050)$\\
$z_6$ & $(-.86466779,-.42329589)$\\
$z_7$ & $(.55580940,-.78606979)$\\
$z_8$ & $(1.00000000,0.00000000)$\\
$z_9$ & $(-.87754040,.47950271)$\\
\end{tabular}\end{minipage}}
\\
\hline\hline
\end{longtable}

%% file: nbody_10.tex
\begin{longtable}{|ll|}
\hline
\multicolumn{2}{|c|}{\emph{Some central configurations for 10 bodies}}\\ 
\hline
\hline \begin{minipage}{0.4\textwidth}
\begin{center} 

\includegraphics[width=\textwidth]{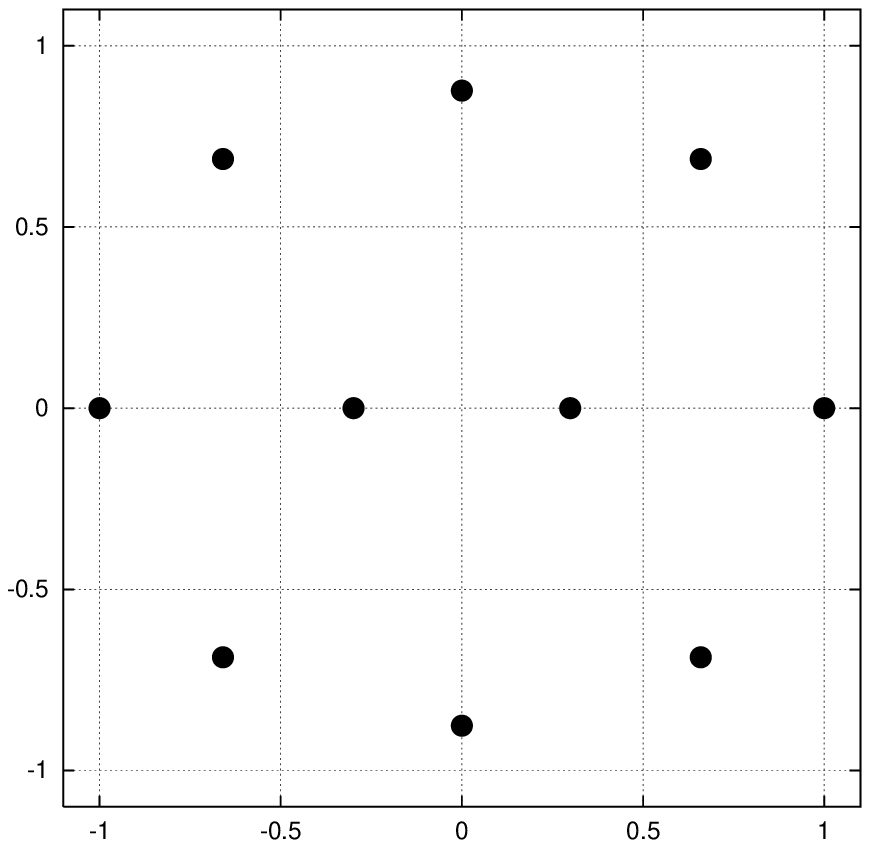} 
\end{center}\end{minipage}
&
\fbox{\begin{minipage}[t]{0.6\textwidth}
\begin{tabular}{r@{ $ = $ }l}
$U\sqrt{I}$ &    $115.60682621 $ \\  
Crit ind & $0$ \\ 
FP ind & $1$ \\ 
Isotropy & $2$\\ 
$z_1$ & $(.29890007,-.00000000)$\\
$z_2$ & $(-.29890007,.00000000)$\\
$z_3$ & $(.00000000,.87611264)$\\
$z_4$ & $(-.00000000,-.87611264)$\\
$z_5$ & $(-.65943000,-.68735185)$\\
$z_6$ & $(.65943000,.68735185)$\\
$z_7$ & $(-.65943000,.68735185)$\\
$z_8$ & $(.65943000,-.68735185)$\\
$z_9$ & $(1.00000000,0.00000000)$\\
$z_{10}$ & $(-1.00000000,-.00000000)$\\
\end{tabular}\end{minipage}}
\\ 
\hline \begin{minipage}{0.4\textwidth}
\begin{center} 

\includegraphics[width=\textwidth]{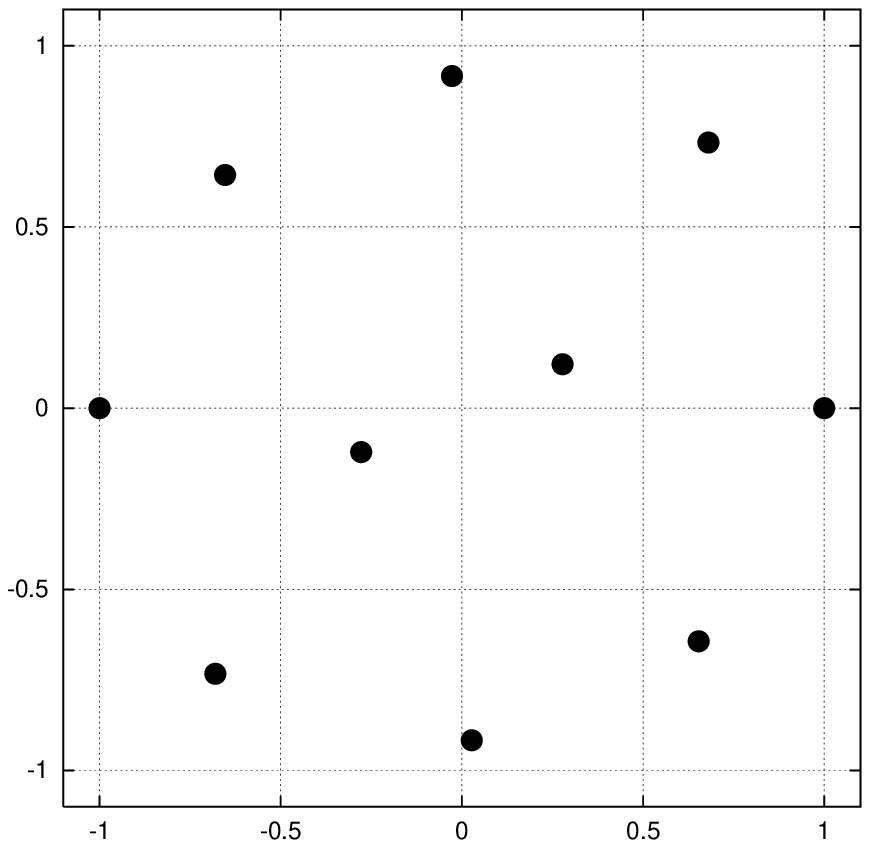} 
\end{center}\end{minipage}
&
\fbox{\begin{minipage}[t]{0.6\textwidth}
\begin{tabular}{r@{ $ = $ }l}
$U\sqrt{I}$ &    $115.60728898 $ \\  
Crit ind & $1$ \\ 
FP ind & $-1$ \\ 
Isotropy & $2$\\ 
$z_1$ & $(.27765563,.12111754)$\\
$z_2$ & $(-.27765563,-.12111754)$\\
$z_3$ & $(-.02716025,.91660090)$\\
$z_4$ & $(.02716025,-.91660090)$\\
$z_5$ & $(.65335522,-.64344530)$\\
$z_6$ & $(-.65335522,.64344530)$\\
$z_7$ & $(-.68027200,-.73295976)$\\
$z_8$ & $(.68027200,.73295976)$\\
$z_9$ & $(1.00000000,0.00000000)$\\
$z_{10}$ & $(-1.00000000,-.00000000)$\\
\end{tabular}\end{minipage}}
\\
\hline\hline
\end{longtable}

%% file: nb.bbl
\begin{thebibliography}{AC99}

\bibitem {A95}
Albouy, Alain.
\emph{Sym\'etrie des configurations centrales
   de quatre corps.} (French) [Symmetry of central configurations of four
   bodies] C. R. Acad. Sci. Paris S\'er. I Math. 320 (1995), no. 2,     
   217--220.

\bibitem{AC98}
Albouy, Alain; Chenciner, Alain.
\emph{Le probl\`eme des $n$
   corps et les distances mutuelles.} (French) [The $n$-body problem and
   mutual distances] Invent. Math. 131 (1998), no. 1, 151--184.

\bibitem{Br71}
Brown, Robert F. \emph{The Lefschetz fixed point theorem.}
   Scott, Foresman and Co., Glenview, Ill.-London 1971 vi+186 pp.

\bibitem{Bu90}
Buck, Gregory. 
\emph{On Clustering in Central Configurations.}
Proceedings of the A.M.S., Vol. 108 n. 3 (1990) 801--810.


\bibitem{Ko01}
Kotsireas, Ilias S. 
\emph{Central configurations in the Newtonian $N$-body problem of celestial mechanics.}
Symbolic computation: solving
equations in algebra, geometry, and engineering (South Hadley, MA, 2000), 
71--98, Contemp. Math., 286, Amer. Math. Soc., Providence, RI, 2001. 


\bibitem{Li01},
Lindstrom, Peter W. 
\emph{The number of planar central configurations is finite when $n-1$ mass positions are fixed.} Trans. Amer. Math.
Soc. 353 (2001), no. 1, 291--311.

\bibitem{MC96}
McCord, Christopher K.
   \emph{Planar central configuration estimates in the $n$-body problem.}
   (English. English summary)
   Ergodic Theory Dynam. Systems 16 (1996), no. 5, 1059--1070.


\bibitem{Mo90}
 Moeckel, Richard.
 \emph{On central configurations.} 
 Math. Z.
   205 (1990), no. 4, 499--517.

\bibitem{Mo01}
Moeckel, Richard. 
\emph{Generic finiteness for Dziobek
   configurations.}
   Trans. Amer. Math. Soc. 353 (2001), no. 11, 4673--4686

\bibitem{moulton}
Moulton, Forest Ray.
\emph{The straight line solutions of the problem of $N$ bodies.}
Ann. Math. II Ser. {12} (1910), 1--17.


\bibitem{PSSY}
Perez-Chavela, Ernesto; Saari, Donald G.; Susin,
Antoni; Yan, Zhiming.
   \emph{Central configurations in the charged three body problem. 
   Hamiltonian dynamics and celestial mechanics (Seattle, WA, 1995)},
   137--153,
   Contemp. Math., 198,
   Amer. Math. Soc., Providence, RI, 1996.

\bibitem{shub}
Shub, m.  
\emph{Appendix to Smale's paper: "Diagonals and relative equilibria"}.
1971
   Manifolds -- Amsterdam 1970 (Proc. Nuffic Summer School) pp. 199--201
   Lecture Notes in Mathematics, Vol. 197.


\bibitem{Sm70I}
 Smale, Steve. 
 \emph{Topology and mechanics. I.} Invent. Math. 10
   (1970), 305--331.

\bibitem{Sm70II}
Smale, Steve. 
\emph{Topology and mechanics. II. The planar
   $n$-body problem.} Invent. Math. 11 (1970), 45--64.
   

   
\bibitem{Sz75}
Szeg\"o, G. \emph{Orthogonal Polynomials},
Fourth Edition, Amer. Math. Soc., Providence, 1975.



\end{thebibliography}
